\documentclass[a4paper,11pt]{amsart}
\usepackage[francais]{babel}
\usepackage[applemac]{inputenc}
\usepackage{epsfig,xypic}
\usepackage{amsthm,amssymb,bbm,a4wide}
\author{Ivan Marin}
\address{Institut de Mathématiques de Jussieu, 175 rue du Chevaleret, F-75013 Paris}
\email{marin@math.jussieu.fr}
\urladdr{http://www.math.jussieu.fr/~marin}
\title{Sur les repr\'esentations de Krammer g\'en\'eriques}
\date{15 juin 2006}

\newtheorem*{theo6}{Th\'eor\`eme 6}
\newtheorem*{theo7}{Th\'eor\`eme 7}
\newtheorem*{theo5a}{Th\'eor\`eme 5a}
\newtheorem*{theo5b}{Th\'eor\`eme 5b}
\newtheorem*{theoA}{Théorème A}
\newtheorem*{theoB}{Théorème B}
\newtheorem*{theoC}{Théorème C}
\newtheorem{theoint}{Th\'eor\`eme}
\newtheorem{lemme}{Lemme}
\newtheorem{definition}{D\'efinition}
\newtheorem{proposition}{Proposition}
\newtheorem*{ccor}{Corollaire}
\newtheorem{cor}{Corollaire}

\def\PPhi{\widetilde{\Phi}}
\def\rg{\mathrm{rg}}
\def\v{\mathbf{v}}
\newcommand{\la}{\lambda}

\def\End{\mathrm{End}}
\def\Ker{\mathrm{Ker}}

\def\SN{\ensuremath{\mathfrak{S}_n}}
\def\C{\ensuremath{\mathbbm{C}}}
\def\Q{\mathbbm{Q}}
\def\Z{\mathbbm{Z}}
\def\R{\mathbbm{R}}
\def\F{\mathbbm{F}}
\def\N{\mathbbm{N}}
\def\gl{\mathfrak{gl}}
\def\sl{\mathfrak{sl}}

\def\k{\mathbbm{k}}
\def\K{\mathbbm{K}}
\def\LL{\mathbbm{L}}

\def\om{\omega}

\def\eps{\epsilon}
\def\g{\mathfrak{g}}
\def\aa{\mathfrak{a}}

\def\into{\hookrightarrow}
\def\ii{\mathrm{i}}

\def\Id{\mathrm{Id}}
\def\lmm{ \left( \begin{array}{cc}}
\def\rmm{\end{array} \right)}

\def\RR{\mathcal{R}}
\def\SS{\mathcal{S}}
\def\Res{\mathrm{Res}}

\def\ss{\sigma}
\def\sss{\tilde{\sigma}}
\def\ee{\tilde{e}}

\begin{document}

\maketitle

\noindent {\bf Abstract.}
We define a representation of the Artin groups of type
ADE by monodromy of generalized KZ-systems which is shown to be
isomorphic to the generalized Krammer
representation originally defined by A.M. Cohen and D. Wales, and independantly
by F. Digne. It follows
that all pure Artin groups of spherical type are residually torsion-free nilpotent,
hence (bi-)orderable.
Using that construction we show that these irreducible representations
are Zariski-dense in the
corresponding general linear group. It follows that all irreducible
Artin groups of spherical type can be embedded as Zariski-dense subgroups
of some general linear group.
As group-theoretical applications we prove properties of non-decomposition
in direct products for several subgroups of Artin groups,
and a generalization in arbitrary types of a celebrated property
of D. Long for the braid groups. We also determine the Frattini and
Fitting subgroups and
discuss unitarity properties of the representations.

\bigskip


\section{Introduction}

Les groupes d'Artin, ou d'Artin-Tits, sont les généralisations
des groupes de tresses associés à chacun des groupes de Coxeter
finis. Ils apparaissent comme espaces d'Eilenberg-MacLane
de variétés algébriques complexes, quotients d'un complémentaire d'arrangement
d'hyperplans (\og espace de configuration \fg) par le groupe de
Coxeter correspondant. Cela permet
notamment de construire une grande partie de leurs représentations par
monodromie de systèmes différentiels sur cet espace de configuration,
par la déformation de représentations du groupe de Coxeter associé.
Il s'agit d'une généralisation des systèmes de Knizhnik-Zamolodchikov
pour le groupe de tresses classique. Les représentations
des groupes d'Artin qui apparaissent de cette façon sont plus faciles
à comprendre que les autres quand les
questions qui se posent peuvent s'exprimer à partir du
système différentiel dont elles proviennent.

Parmi les représentations du groupe de tresses classique,
une représentation retient particulièrement l'attention depuis
quelques années. Il
s'agit d'une représentation irréductible de l'algèbre de
Birman-Wenzl-Murakami,
introduite en 1989, qui a été particulièrement étudiée par D. Krammer
et dont il a montré en 2000,
en concurrence avec S. Bigelow, qu'elle était fidèle (cf. \cite{KRAM,BIG}).
Cette
représentation apparaît également dans les travaux de R. Lawrence \cite{LAW} en 1991,
c'est pourquoi elle est parfois appelée représentation de Lawrence-Krammer.
On montre facilement que ces différentes constructions donnent lieu
à des représentations équivalentes. Pour des comparaisons explicites, on
pourra consulter \cite{PP,ZINNO}.

Les représentations de l'algèbre de Birman-Wenzl-Murakami
font partie des représentations qui peuvent s'obtenir par monodromie.
D'autre part, la représentation de Krammer et la preuve de
sa fidélité a été généralisée, par Cohen et Wales \cite{CW} en types ADE,
et par Digne \cite{DIGNE} pour tous les types cristallographiques. Comme il s'agit essentiellement
de la seule représentation de ces groupes dont la fidélité soit
connue, son étude présente un intérêt particulier.

Le but de ce travail est, d'une part, de construire un
système différentiel sur l'espace de configuration en types ADE,
c'est-à-dire une \og représentation infinitésimale \fg\ (prop. \ref{proprepinf}), dont
la monodromie est la représentation de Krammer en types ADE
(théorème A).

Cette construction nous
permet de définir une forme quadratique infinitésimale
qui explique naturellement l'apparition de structures unitaires sur la représentation
de Krammer en type A, et conjecturalement en types D et E
(cf. \cite{REPTHEO,DIEDRAUX}).

D'autre part, elle permet de déterminer l'enveloppe algébrique des
représentations de Krammer, et de déterminer la décomposition
de ses puissances tensorielles. Nous montrons ici (théorème B) que le
groupe d'Artin considéré est Zariski-dense dans le groupe linéaire
associé à la représentation de Krammer, et
qu'il en est déjà de même pour certains de ses sous-groupes.

La fidélité de cette représentation permet d'obtenir des résultats
nouveaux de théorie des
groupes, par exemple la résiduelle nilpotence des
groupes d'Artin purs.
De ce résultat de densité pour la topologie de Zariski
nous déduisons de plus des conséquences sur les possibles décompositions en produit direct
de sous-groupes, notamment des sous-groupes d'indice fini, des groupes d'Artin. Nous en
déduisons enfin les sous-groupes de Fitting et de Frattini, et généralisons
en tous types une propriété, due à D. Long, des groupes de tresses.

\bigskip

La section suivante présente les notations utilisées et détaille,
sans les démontrer, les principaux résultats issus de
ce travail.
En section 3 nous définissons l'action infinitésimale,
étudions ses propriétés et introduisons une forme quadratique
remarquable $(\ | \ )$. En section 4 nous démontrons que la monodromie
de cette action est isomorphe à la représentation de Krammer généralisée,
en section 5 nous démontrons sa Zariski-densité et présentons
des applications immédiates. La section 6 développe les conséquences
de ce résultats pour les groupes d'Artin-Tits sphériques
irréductibles de type arbitraire. Enfin, les sections 7 et 8 contiennent une
étude des différents cas $A_n$ et $D_n$, qui sert notamment
à déterminer le discriminant de $(\ | \ )$.

\bigskip

\noindent {\bf Remerciements.} Les calculs ont bénéficié des
ressources du centre de calcul Médicis de l'Ecole Polytechnique. Je remercie
de nombreuses personnes pour des échanges
qui m'ont stimulé dans ce travail, notamment
Y. de Cornulier, C. Cornut, F. Courtes, J. Crisp, F. Digne, D. Krammer, L. Paris.
Je remercie en particulier Y. de Cornulier pour de précieux commentaires
sur une version préliminaire, ainsi que F. Digne pour avoir rectifié
les formules de \cite{DIGNE}.

\section{Notations et présentation des résultats}

\subsection{Groupes d'Artin et représentation de Krammer}

On note $W$ un groupe de Coxeter fini irréductible, considéré comme
groupe de réflexion d'un $\R$-espace vectoriel de dimension  $n$,
où $n$ désigne le rang de $W$. On dira que $W$ est de type
$ADE$ s'il est d'un des types $A_n$, $D_n$
ou $E_n$, avec $n \in \{ 6,7,8 \}$. Rappelons que cette
condition est équivalente à la propriété que
toutes les réflexions de $W$ se trouvent dans
la même classe de conjugaison. D'autre part, rappelons
que le groupe de Coxeter de type $I_2(m)$ est le groupe
diédral d'ordre $2m$, et que celui de type $A_{n-1}$
est le groupe symétrique $\mathfrak{S}_n$.

On notera $\RR$ l'ensemble des réflexions de $W$. A
chacune de ces réflexions $s$ est associé un hyperplan
$H = \Ker(s-1)$ de l'espace vectoriel $\R^n$ considéré,
donc à $W$ est associé un arrangement d'hyperplans réels.
Complexifiant $\R^n$ en $\C^n$, on en déduit
une action de $W$ sur $\C^n$ et un arrangement
$\mathcal{A}$ d'hyperplans complexes. On note $X_W$
le complémentaire dans $\C^n$ de cet arrangement d'hyperplans.

Les groupes d'Artin-Tits dits \og sphériques \fg\ sont
les groupes fondamentaux $B = \pi_1(X_W/W)$ associés
à de tels groupes de Coxeter finis. Ils
sont dits irréductibles lorsque le groupe de Coxeter
correspondant l'est. Notant $P = \pi(X_W)$
le groupe d'Artin pur correspondant, on a une
suite exacte naturelle
$$
1 \to P \to B \stackrel{\pi}{\to} W \to 1
$$
Lorsque $W$ est de type $A_{n-1}$, le groupe $B$ (resp. $P$)
est le groupe des tresses (pures) à $n+1$ brins.

Supposons $W$ défini par un système $\{ s_1,\dots,s_n \}$ de réflexions
simples. Par référence au
diagramme de Coxeter de $W$, on dira que $s_i$ est relié à $s_j$
si les noeuds $i$ et $j$ du diagramme de Coxeter sont reliés à une
arête : cela signifie que le sous-groupe $<s_i,s_j>$ n'est pas de type $I_2(2)$.
On peut naturellement associer à ce système de réflexions simples
un système d'antécédents $\sigma_1,\dots,\sigma_{n}$
dans $B$,
qui sont appelés des générateurs d'Artin de $B$ : ils
engendrent $B$ et les relations entre ces éléments
se lisent sur le diagramme de Coxeter de $W$. Précisément,
$B$ est algébriquement défini par générateurs $\ss_1, \dots, \ss_n$
et relations $\ss_i \ss_j = \ss_j \ss_i$ si $s_i$ et $s_j$ ne
sont pas reliées, et
$$
\underbrace{\ss_i \ss_j \ss_i \dots }_m
= \underbrace{\ss_j \ss_i \ss_j \dots }_m
$$
si $<s_i,s_j>$ est un groupe diédral de type $I_2(m)$.
Le morphisme surjectif $\pi : B \to W$ de noyau $P$
vérifie donc $\pi(\ss_i) = s_i$ pour tout $i \in [1,n]$.

Introduisons
également
$C_W^k$, pour $k \geq 1$, le sous-groupe de $P$
engendré par $\sigma_1^{2k}, \dots, \sigma_n^{2k}$.
Un résultat de J. Crisp et L. Paris (cf. \cite{CRPA}) montre que $C_W^1$,
donc tout $C_W^k$ pour $k \geq 1$, est un groupe localement libre
en ces générateurs, et ceci quel que soit $W$ : notant
$x_i = \sigma_i^{2k}$, cela signifie que les
seules relations entre ces éléments sont
$x_i x_j = x_j x_i$ si les noeuds du diagramme
de Coxeter de $W$ correspondant à $i$ et à $j$ ne sont
pas reliés entre eux.

Soit $\mathcal{L} = \Q[q,q^{-1},t,t^{-1}]$, et supposons que
$W$ est de type ADE. A la suite
des travaux de Krammer, Cohen et Wales d'une part,
F. Digne de l'autre, ont indépendamment construit
une représentation linéaire fidèle
$ R_{\mathrm{K}} : B \to GL_N(\mathcal{L})$, où $N$ désigne
le nombre de réflexions de $W$, que nous appelerons désormais représentation
de Krammer. Nous ne rappellerons pas les formules
qui définissent cette représentation, d'une part parce qu'elles
sont compliquées et peu éclairantes, d'autre part parce
que nous ne les utiliserons pas ici. Nous renvoyons à \cite{CGW}
pour cette description : nos paramètres sont reliés
aux paramètres $r$ et $l$ de \cite{CGW} par $r = q^{-1}$
et $t = l^{-1}$.

\subsection{Algèbre d'holonomie et représentations de monodromie}

Une manière désormais classique d'obtenir des re\-pré\-sen\-ta\-tions des
groupes d'Artin consiste à étudier des re\-pré\-sen\-ta\-tions
de l'algèbre de Lie d'holonomie de l'arrangement d'hyperplans $\mathcal{A}$,
telle que définie et étudiée par Kohno, notamment dans \cite{KOHNO}. Les
exemples en type A sont légion. En dehors du type A, on pourra par
exemple consulter \cite{CHEREDNIK,BMR,MT}.

Fixant un corps $\k$ de caractéristique 0, cette algèbre de Lie $\g$
définie sur $\k$ est présentée par une famille de générateurs $t_s, s \in \RR$, soumise
aux relations $[t_s,t_Z] = 0$ où
\begin{itemize}
\item $Z \subset \C^n$ parcourt les sous-espaces de
codimension 2 intersections d'éléments de $\mathcal{A}$
\item  $s$ parcourt les éléments de $\RR$ tels que $s_{|Z} = 
\Id_{Z}$
\item $t_Z$ désigne la somme des $t_u$ pour $u_{|Z} = 
\Id_{Z}$.
\end{itemize}
Pour tout $\SS \subset \RR$ on notera $\g_{\SS}$ la sous-algèbre
de Lie de $\g$ engendrée par les éléments $t_s, s \in \SS$.

Notons $\mathcal{M} = \k[[h]]$ l'anneau des séries
formelles en une indéterminée sur $\k$ et $\K = \k((h))$
son corps des fractions. On note alors $\g^h = \g \otimes \mathcal{M}$.

A tout $s \in \R$, on associe une forme linéaire $\alpha_s$
sur $\C^n$ de noyau $\Ker(s-1)$.
La 1-forme
différentielle $\mathrm{d} \log \alpha_s = 
\frac{\mathrm{d} \alpha_s}{\alpha_s}$
ne dépend pas de la forme particulière choisie, et les
relations de $\g$ impliquent que la 1-forme, à valeur dans $\g^h$,
$$
\om = h \sum_{s \in \mathcal{\RR}} t_s \mathrm{d} \log \alpha_s \in \Omega^1(X_W) \otimes \g^h
$$
En conséquence, si $\k = \C$, toute représentation
$\rho : \g \to \gl(V)$, où $V$ est un $\k$-espace vectoriel,
donne lieu par monodromie à une représentation de $P$
sur $V^h = V \otimes \mathcal{M}$. Une telle représentation
ne dépend pas du point-base choisi à isomorphisme près.
De plus, si $V$ est muni d'une action
de $W$ telle que
$$
\forall s \in \RR \ \ \ \forall w \in W \ \ \ w \rho(t_s) w^{-1} = \rho(t_{wsw^{-1}})
$$
on dira que $\rho$ est équivariante. En d'autre termes, elle
s'étend en une représentation de $\mathfrak{B} = \k W \ltimes \mathsf{U} \g$,
où $\mathsf{U} \g$ désigne l'algèbre enveloppante universelle de $\g$.
La représentation de monodromie s'étend alors en une représentation $S : B \to GL(V^h)$.
Les propriétés des représentations construites de cette façon
qui nous seront utiles sont les
suivantes :
\begin{enumerate}
\item Si $\ss$ est un générateur d'Artin et $s = \pi(\ss)$, alors $S(\ss)$ est conjugué à $s \exp( \ii \pi h \rho(t_s))$.
\item Si $\rho$ est (absolument) irréductible, alors $S$ est (absolument) irréductible
\end{enumerate}
La première de ces propriétés est classique ; pour la deuxième on
pourra consulter \cite{REPTHEO}.

\subsection{Résultats sur la représentation de Krammer}

En section 3 nous construisons, pour tout groupe de
Coxeter $W$ irréductible fini de type $ADE$ et tout $m \in \k$,
une représentation de $\g$ dans le sous-espace vectoriel
de base $\RR$, que nous appelons représentation
de Krammer infinitésimale. Cette représentation est compatible
avec l'action naturelle par conjugaison de $W$ sur $\RR$,
et permet donc de construire par monodromie, si $\k = \C$, une
représentation $R : B \to GL_N(\mathcal{M}) \subset GL_N(K)$.

Un plongement de $\mathcal{L}$ dans $\mathcal{M} = \k[[h]]$, donc
dans $\K= \k((h))$, est défini par la donnée des images de $q$ et $t$,
qui doivent être inversibles et algébriquement indépendantes.
C'est le cas notamment si $q = e^h$ et $t = e^{m h}$
pour $m \not\in \Q$.

Les formules définissant cette \og représentation de Krammer
infinitésimale \fg\ sont nettement plus simples que
les descriptions connues de la représentation
de Krammer dans le cas général. Nous montrons en section 4
(théorème A) que, si $N$
désigne le nombre de réflexions de $W$,

\begin{theoint} Pour tout $W$ de type $ADE$
et $m \not\in \Q$,
la représentation de Krammer $R_{\mathrm{K}}$ plongée dans $GL_N(\K)$
par $q = e^h$ et $q= e^{m h}$ est isomorphe à $R$.
\end{theoint}

Nous montrons en section 3 que $\rho$, donc $R$, est irréductible.
Nous montrons également que $\rho : \g \to \gl_N(\k)$ est surjective,
et même que $\rho(\g_{\SS}) = \gl_N(\k)$ pour $\SS$
composé des racines simples de $W$. Une conséquence (section 5.1,
théorème B) est que

\begin{theoint} Pour tout $W$ de type $ADE$ et $m \not\in \Q$, la représentation
de Krammer $R_{\mathrm{K}}$ plongée dans $GL_N(\K)$
par $q = e^h$ et $q= e^{m h}$ est d'image Zariski-dense.
\end{theoint}

\subsection{Conséquences sur les groupes d'Artin-Tits}

Une conséquence du premier théorème est le fait que
$P$ se plonge, pour tout $W$ de type $ADE$, dans
$$
GL_N^0(\K) = 1 + h M_N(\mathcal{M}) \subset GL_N(\K)
$$
qui est un groupe résiduellement nilpotent-sans-torsion.
Comme tous les groupes d'Artin-Tits sphériques
peuvent se plonger dans un groupe de type $ADE$ d'après \cite{CRISP},
on en déduit immédiatement (corollaire du théorème A)
\begin{theoint}
Tout groupe d'Artin-Tits sphérique pur est résiduellement
nilpotent-sans-torsion, et en particulier biordonnable.
\end{theoint}
La résiduelle nilpotence-sans-torsion avait déjà été obtenue en types
$A,B$ (et $I_2(m)$) par Falk et Randell \cite{FALK} et ultérieurement
en type $D$ par Leibman et Markushevich \cite{LM}.
Remarquons qu'en posant artificiellement
$q = e^h$ et $t = e^{\sqrt{2} h}$ dans les formules de la
représentation de Krammer, il est possible d'obtenir directement, mais de
façon un peu magique, ce
résultat. Nous renvoyons à \cite{resnil} où nous avons utilisé cette méthode.
Le fait que tout groupe ré\-si\-duel\-le\-ment nilpotent-sans-torsion est biordonnable
est démontré par exemple dans \cite{PASS}.

Une con\-sé\-quen\-ce du deuxième théo\-rè\-me est que tout groupe
d'Artin-Tits sphé\-ri\-que
ir\-ré\-duc\-ti\-ble de type $ADE$ est sous-groupe Zariski-dense d'un groupe
linéaire. Il en est donc de même pour les groupes
de type $B_n$, qui s'identifient à des sous-groupes d'indice fini
de groupes de type $A_n$. Un résultat analogue s'obtient facilement pour
les groupes diédraux, qui se plongent dans $GL_2(K)$ d'après un
résultat de G. Lehrer et N. Xi \cite{LEHRER}. Le cas des groupes exceptionnels
$F_4,H_3,H_4$ nécessite un peu plus de travail, pour compléter la
démonstration (section 6) du

\begin{theoint}
Pour tout groupe d'Artin-Tits sphérique irréductible $G$, il
existe $N \in \N$ et un corps $K$ tel que
\begin{enumerate}
\item $B$ s'identifie
à un sous-groupe Zariski-dense de $GL_N(K)$
\item Les sous-groupes $C_W^k$ sont également Zariski-denses.
\end{enumerate}
\end{theoint}

Soulignons que ce résultat n'était pas connu pour
le groupe de tresses \og classique \fg.

\subsection{Conséquences de la Zariski-densité}

Le théorème 4 a plusieurs conséquences remarquables de type \og théorie des groupes \fg\ 
sur les groupes d'Artin-Tits sphériques irréductibles. Rappelons tout
d'abord que, si $B$ est un tel groupe, son centre $Z$ est cyclique infini.

En premier lieu, cela permet de retrouver
facilement un résultat de L. Paris.

\begin{theo5a} (L. Paris, \cite{PARIS}) Soit $B$ un groupe d'Artin-Tits irréductible,
et $Z$ son centre.
Si $B = U \times V$ est
une décomposition de $B$ en produit direct, alors $U \subset Z$
ou $V \subset Z$.
\end{theo5a}

De plus, cette nouvelle démonstration (section 6.2) montre que de nombreux autres groupes
ont cette propriété, par exemple les sous-groupes d'indice fini de
$B$.

\begin{theo5b} Soit $B$ un groupe d'Artin-Tits sphérique irréductible, et $G$ un
sous-groupe d'indice fini de $B$. Si $G = U \times V$ est
une décomposition de $G$ en produit direct, alors $U \subset Z$
ou $V \subset Z$.
\end{theo5b}

Plus généralement, on en déduit une généralisation
à tous les groupes d'Artin-Tits d'un théorème
fameux de D. Long sur le groupe de tresses. Dans \cite{LONG},
il est en effet démontré à l'aide de la classification
de Nielsen-Thurston des difféomorphismes que les sous-groupes
sous-normaux $Q$ des groupes $G$ qui sont des mapping class groups
ou des groupes de tresses vérifient la propriété suivante,
où l'on a noté $Z$ le centre de $G$.
\begin{center}
(L) Si $H_1,H_2$ sont deux sous-groupes distingués propres de
$Q$ non inclus dans $Z$, alors $H_1 \cap H_2$ n'est pas inclus
dans $Z$. 
\end{center}
Cette propriété a des conséquences importantes en théorie des
représentations, déjà pour $Q = G$. Elle implique notamment
que 
\begin{enumerate}
\item on ne peut essentiellement pas construire de représentation linéaire fidèle
par somme directe de représentations non fidèles
\item pour vérifier qu'une représentation d'un groupe de tresses
est fidèle, il suffit de vérifier la fidélité sur le centre
et le sous-groupe libre naturel.
\end{enumerate}

Une conséquence du théorème 4 est que cette propriété se généralise
à tous les sous-groupes d'indice fini des groupes d'Artin-Tits sphériques

\begin{theo6} Soit $B$ un groupe d'Artin-Tits sphérique irréductible,
et $G$ un sous-groupe d'indice fini de $B$. Alors tous les sous-groupes
sous-normaux de $G$ vérifient (L).
\end{theo6}

Plus précisément, nous étudions en section 6.3 les propriétés
des groupes sans torsion, appelés ici \og fortement linéaires \fg, qui se
plongent dans un groupe linéaire $GL_N(K)$, pour un certain corps $K$
et un certain entier $N$, de telle façon que leur
adhérence pour la topologie de Zariski contienne $SL_N(K)$. Les groupes d'Artin-Tits
sphériques irréductibles en font partie d'après le théorème 4,
ainsi que leurs sous-groupes d'indice fini, soit par des arguments
généraux (prop. \ref{flinfini}), soit parce qu'ils contiennent les sous-groupes
$C_W^k$ pour $k$ assez grand. Le théorème 6
découle alors d'une propriété générale de ces groupes fortement
linéaires (théorème C, section 6.3).

En particulier, ce théorème montre que, pour déterminer la
fidélité d'une représentation d'un groupe d'Artin-Tits
sphérique irréductible $B$, il suffit de vérifier la
fidélité de sa restriction au centre et à n'importe
lequel de ses sous-groupes sous-normaux. Or, tous ces
groupes sauf peut-être un nombre fini d'entre eux
(correspondant aux types exceptionnels) admettent
un sous-groupe sous-normal qui est libre, donc auquel
des méthodes de ping-pong sont susceptibles de s'appliquer.

En effet, on a un morphisme surjectif $\pi : B \to W$
vers le groupe de Coxeter correspondant, et
$P = \Ker \pi$ s'identifie à $\pi_1(X)$, pour
$X$ le complémentaire d'un certain arrangement d'hyperplans.
En types $A,B$ et $I_2(m)$ cet arrangement est super-résoluble
ou encore \og fiber-type \fg, suivant la terminologie de \cite{ORLIKTERAO},
donc il existe une fibration $X \to Y$, de fibre $\C \setminus \{ m \mbox{ points} \}$
et dont l'espace de base est de $\pi_2$ nul. Cela
induit un plongement $\pi_1( \C \setminus \{ m \mbox{ points} \}) \into
\pi_1(X)$ comme sous-groupe normal, donc comme sous-groupe sous-normal de $B$.
En type $D$ il existe également une telle fibration, due à Brieskorn,
introduite dans \cite{BRIESK} et particulièrement étudiée par Markusevitch
et Leibman dans \cite{LM}. Les seuls cas dans lesquels l'existence
d'un tel sous-groupe sous-normal libre reste à notre connaissance ouverte sont donc les
types exceptionnels $H_3,H_4,F_4,E_6,E_7,E_8$. 

Enfin, ces résultats permettent d'obtenir facilement les sous-groupes
de Fitting $\mathrm{Fit}(G)$ et de Frattini $\Phi(G)$ de ces groupes.
En type $A$, le groupe $\Phi(B)$ a été déterminé par Long dans \cite{LONG}.

\begin{theo7}
Soit $Q$ un sous-groupe sous-normal d'un groupe d'Artin-Tits sphérique
ir\-ré\-duc\-ti\-ble $B$. Alors $\mathrm{Fit}(Q) = Z(Q)$ et,
si $Q$ est finiment engendré, alors $\Phi(Q) = \{ e \}$.
En particulier, $\Phi(B) = \{ e \}$. 
\end{theo7}
\begin{proof} Le fait que $\mathrm{Fit}(Q) = Z(Q)$ découle du
théorème C et de son deuxième corollaire. Supposons $Q$ finiment
engendré. D'après le corollaire 3 du théorème C on a $\Phi(Q)
\subset Z(Q) \subset Z(B)$. Soit $\ell : B \to \Z$ le morphisme
longueur associé aux générateurs d'Artin, et $\check{\ell} : Q \to \Z$
sa restriction à $Q$. Si $\check{\ell}(Q) = \{ 0 \}$, alors
$\Phi(Q) \subset \Ker \ell$ d'où $\Phi(Q) = \{ e \}$ car
$Z(B)$ est engendré par un élément de longueur non nulle. Sinon
$\check{\ell}(Q) = m \Z$ avec $m \in \N^*$. Or, pour tout nombre premier
$p$ ne divisant pas $m$, le sous-groupe $\check{\ell}^{-1}(pm\Z)$
est maximal dans $Q$ donc $\Phi(Q) \subset \check{\ell}^{-1}(p m \Z)$.
Comme $m \Z$ est monogène $\check{\ell}^{-1}(pm\Z)$ est maximal
dans $Q$ donc $\Phi(Q) \subset \Ker \ell$ et on conclut comme
précédemment.
\end{proof}

\subsection{Discussions sur l'unitarité.}

Un dernier aspect de cette approche infinitésimale de la représentation
de Krammer concerne la présence d'une forme quadratique, simple
à définir et naturelle, qui est \og invariante \fg\ par l'action
infinitésimale en type $ADE$. De telles formes conduisent à l'existence
de formes unitaires invariantes sur la représentation
de Krammer générique par des méthodes de monodromie
algébrique, développées dans \cite{DRINFELD,ENRIQUEZ,DIEDRAUX} en types $A$, $B$, et $I_2(m)$.
De telles méthodes n'ont pas encore été développées en type $D$ et $E$, donc
cela ne permet pour l'instant que de démontrer
l'unitarisabilité {\it a priori} (prop. \ref{propunita}) de
la représentation de Krammer en type $A$, explicitée par R. Budney dans \cite{BUDNEY} par de
toutes autres méthodes. Néanmoins, il nous a semblé utile de décrire explicitement
cette forme quadratique, très vraisemblablement à l'origine de formes
unitaires invariantes, pour les différents types. Notamment, nous
déterminons les facteurs irréductibles de son discriminant.

\section{Représentations de Krammer infinitésimales}

\subsection{Définitions des opérateurs $\tau_s$}

On note $a \bullet b = a b a^{-1}$ et on fixe un corps $\k$ de caractéristique 0.
Soit $W$ un groupe de Coxeter irréductible fini de type $ADE$ de rang $n$,
identifié à un sous-groupe de réflexions de $\C^n$, $\RR$ l'ensemble de ses réflexions
et
$\g$ l'algèbre de Lie d'holonomie associée.

Soit $V$ un $\k$-espace vectoriel de base $\{ v_s, s \in \RR\}$,
muni de l'action de $W$ définie par $w.v_s = v_{ w \bullet s}$.
Le vecteur $\v$ égal à la somme des $v_s$ pour $s \in \RR$ est $W$-invariant,
et engendre $V^W$ parce que toutes les réflexions sont conjuguées dans $W$.

On va définir, pour tout $m \in \k$, une famille d'éléments $\tau_s \in \End(V)$
pour $s \in \RR$, dont on montrera ensuite qu'ils
permettent de définir une représentation équivariante de $\g$.
\begin{definition} Pour tout $s \in \RR$ on note $\tau_s \in \End(V)$
l'endomorphisme défini par
$$
\tau_s . v_u = \left\lbrace \begin{array}{lcl}
v_{s \bullet u } - v_s & \mbox{ si } & su \neq us \\
v_{ u }  & \mbox{ si } & su = us, u \neq s \\
m v_s & \mbox{ si } & u=s  \\
\end{array} \right.
$$
\end{definition}
On identifie $w \in W$ à son action naturelle
sur $V$.
\begin{lemme} Pour tous $s \in \mathcal{R}$, $w \in W$, on a
$ w \tau_s w^{-1} = \tau_{w \bullet s}$
\end{lemme}
\begin{proof}
On note $\delta_{u,s}$ le symbole de Kronecker, i.e.
$\delta_{u,s} = 1$ si $u=s$, $0$ sinon, et
$\delta'_{u,s} = 1$ si $u \neq s$ et $us = su$, $0$ sinon.
Soit $u \in \mathcal{R}$.
On remarque $
\tau_s . v_u = v_{s \bullet u} - v_s + m\delta_{u,s} v_s + \delta'_{u,s} v_s,
$. On a
$$
\begin{array}{lcl}
 w \tau_s w^{-1}. v_u & = &  w \tau_s .v_{w^{-1} \bullet u} \\
& = & w( v_{s w^{-1} \bullet u} - v_s + m \delta_{s,w^{-1} \bullet u}
v_s + \delta'_{s, w^{-1} \bullet u} v_s )\\
& = & v_{ws w^{-1} \bullet u} - v_{w \bullet s} + m \delta_{s,w^{-1} \bullet u}
v_{w \bullet s} + \delta'_{s, w^{-1} \bullet u} v_{w \bullet s} \\
& = & v_{(w \bullet s) \bullet u} - v_{w \bullet s} + m \delta_{w \bullet s,u}
v_{w \bullet s} + \delta'_{w \bullet s,  u} v_{w \bullet s} \\
& = & \tau_{w \bullet s} . v_u 
\end{array}
$$
\end{proof}

Remarquons que, identifiant toujours $s \in \RR$ à son action naturelle
sur $V$, on a $\tau_s = s - p_s$ avec
$$
p_s . v_u = \left\lbrace \begin{array}{lcl}
 v_s & \mbox{ si } & su \neq us \\
0  & \mbox{ si } & su = us, u \neq s \\
(1-m) v_s & \mbox{ si } & u=s  \\
\end{array} \right.
$$
On en déduit $w p_s w^{-1} = p_{w \bullet s}$.

En particulier, $s$ et $p_s$ commutent et on a $p_s^2 = (1-m) p_s$.
On vérifie facilement $s p_s = p_s s = p_s$ pour tout $s \in \RR$
donc $p_s \tau_s = \tau_s p_s = m p_s$. On en déduit
$$
\tau_s^3 - \tau_s -m \tau_s^2 + m = 0,\ \ \ s = \frac{-1}{m+1} \tau_s^2 + \tau_s +
\frac{1}{m+1} \mbox{ si } m \neq -1
$$
et donc les spectres comparés de $\tau_s$ et $s$ sont
$$
\begin{array}{|c||c|c|c|}
\hline
 \tau_s & 1 & -1 & m \\
\hline
s & 1 & -1 & 1 \\
\hline
\end{array}
$$
\subsection{Une forme quadratique invariante}

Introduisons la forme bilinéaire symétrique telle que $(v_r|v_s) = m-1$ si
$r=s$, $(v_r|v_s) = -1$ si $rs \neq sr$ et $(v_r|v_s)=0$ sinon.
Il est clair qu'elle est invariante pour l'action de $W$.

\begin{proposition} \label{fquad} Lorsque $(\ | \ )$ est non dégénérée
et $m \neq 1$, pour tout $s \in \RR$,
$\frac{1}{1-m} p_s$ est le projecteur orthogonal sur
$\k v_s$ par rapport à cette forme et $\tau_s$ est autoadjoint.
Enfin, $(\ | \ )$ est non dégénérée pour des valeurs génériques de $m$
et, si $\k \subset \R$, c'est un produit scalaire pour $m$ assez grand.
\end{proposition}
\begin{proof}
Comme chaque réflexion
est à la fois orthogonale par rapport à $(\ | \ )$ et involutive, on en
déduit
que $r$ est autodajoint pour tout $r \in \RR$.
Comme on a $p_r^2 = (1-m) p_r$ et
$\mathrm{Im} p_r = \k v_r$, il
suffit de montrer que $p_r$ est autoadjoint pour avoir la
première partie de l'énoncé.
On a $p_r (v_s) = -(v_r|v_s)v_r$ d'où $(p_r(v_s)|v_{s'}) = -(v_r|v_s)(v_r|
v_{s'}) = (v_s | p_r(v_{s'}))$ et
$p_r$ est bien autoadjoint. Le déterminant de la matrice définissant
$(\ | \ )$ dans la base naturelle est un polynôme unitaire non constant
en $m$, donc cette
forme quadratique est bien non dégénérée pour $m$ générique, et si
$\k \subset \R$ ses mineurs principaux sont également des
polynômes unitaires non constant en $m$, donc sont strictement positifs
si $m$ est suffisamment grand. On en déduit d'après
le critère de Sylvester que, pour $m$ suffisamment
grand, $(\ | \ )$ est un produit scalaire.
\end{proof}

Quand $m$ n'est pas générique, on peut dire a priori la chose suivante.
Par abus de langage, on appellera discriminant de $(\ | \ )$ le
déterminant de sa matrice dans la base formée des $v_s, s \in \RR$
pour un ordre donné sur $\RR$. Il s'agit donc d'un polynôme unitaire en $m$
à coefficients rationnels de degré $\# \RR$.

Pour tout $s \in \RR$, introduisons les nombres
$$
\begin{array}{lcl}
c(W) & = & \# \{ s' \in \mathcal{R}\setminus \{ s \}\   | \ s's = ss' \}\\
c'(W) & = & \# \{ s' \in \mathcal{R} \ | \ s's \neq ss' \} 
\end{array}
$$
qui ne dépendent pas du choix de $s$ puisque toutes les
réflexions sont conjuguées entre elles, et notons $V_1 \oplus \dots \oplus V_r$ une décomposition en
irréductibles de $V$ en tant que $W$-module, avec
$\dim V_i \leq \dim V_{i+1}$. 
On constate par examen de chacune des séries que cette décomposition est
sans multiplicités, c'est-à-dire $V_i \not\simeq V_j$ pour $i\neq j$.
D'autre part, comme $W$ est de type ADE, toutes ces représentations
sont absolument irréductibles et définies sur $\Q$.
Ainsi $(\ | \ )$ appartient au
carré symétrique de $V^*$ en tant que $W$-module, est $W$-invariant
et sa dépendance en $m$ est affine. 
Comme $V$ est autodual en tant que $W$-module, on peut identifier
canoniquement $(\ | \ )$ à un élément, dépendant affinement de $m$,
de $\bigoplus_{i=1}^r \End_W(V_i) \simeq \k^r$
et son discriminant se décompose comme produit de $r$ facteurs
de la forme $Q_i(m)^{\dim V_i}$, avec $Q_i$ polynôme de degré 1 en $m$.
On peut d'autre part supposer que $V_1$
est la représentation triviale, engendrée par la somme $\v$
des $v_s$ pour $s \in \RR$. On a $(\v | \v) = (\# \RR) (m-\# \RR + c(W))$,
donc $m- \# \RR + c(W)$ apparaît comme facteur irréductible du
discriminant de $(\ | \ )$. Enfin, $V \simeq V^*$ en tant que
$W$-module, donc $V^* \otimes V^*$ s'identifie à la représentation de
permutation naturelle sur $\RR \times \RR$. Comme cette
action comprend au moins trois parties stables symétriques
$\{ (v_s,v_s) \ | \ s \in \RR \}$, 
$\{ (v_s,v_u) \ | \ <s,u> \simeq I_2(3) \}$ et
$\{ (v_s,v_u) \ | \ <s,u> \simeq I_2(2)\}$, toutes non vides si $n \geq 3$,
on en déduit $r \geq 3$ si $n \geq 3$. On peut vérifier au cas par cas
que $r = 2$ pour $n=2$, $r = 3$ pour $W$ de type $A_n$ avec $n \geq 3$
et $W$ de type $E_6,E_7,E_8$, $r=4$ pour $W$ de type $D_n$
avec $n \geq 5$ et $r=5$ pour $W$ de type $D_4$, et déduire
la valeur du discriminant de cette décomposition explicite.
Nous obtenons plus loin la valeur du discriminant pour chacun
de ces types par une méthode alternative, qui permet de plus de décrire la
restriction de $V$ à des sous-algèbres de Lie
remarquables de $\g$. Une fois les facteurs $Q_i$ obtenus, une
condition nécessaire et suffisante de définie positivité pour $\k = \R$
est alors donnée par $Q_i(m) > 0$ pour tout $1 \leq i \leq r$.
Nous verrons type par type que cette condition s'exprime de la façon suivante.
\begin{proposition} Si $\k = \R$, la forme $(\ | \ )$ est définie
positive si et seulement si $m > \# \RR - c(W)$.
\end{proposition}

\begin{table}
$$
\begin{array}{|l|c|l|}
\hline
\mathrm{Type} & \mathrm{Discriminant} & \mbox{Déf. $>0$ ssi} \\
\hline
E_6 & (m-21)(m-3)^{20}(m+3)^{15} & m > 21 \\
\hline
E_7 & (m-33)(m-5)^{27}(m+3)^{35} & m > 33 \\
\hline
E_8 & (m-57)(m-9)^{35}(m+3)^{84} & m > 57 \\
\hline
\end{array}
$$
\caption{Discriminant de $(\ | \ )$ en type E.}
\label{tableE}
\end{table}

Cela découlera des propositions \ref{discrA} et \ref{discrD} pour
les types classiques. Pour les types $E$, cela découle d'un calcul
explicite (cf. table \ref{tableE}).

\subsection{Restriction aux sous-espaces $E_X$}
Soit $X$ une
intersection d'hyperplans de réflexion de codimension 2. On note
$$
\RR_X = \{g \in \RR \ | \ \forall x \in X \ \ g(x) = x \},\ \
W_X = \{g \in W \ | \ \forall x \in X \ \ g(x) = x \}, \ \
\tau_X = \sum_{g \in \RR_X} \tau_g,
$$
et $E_X$ le sous-espace de $V$ engendré par les $v_s$ pour
$s \in \RR_X$. C'est un sous-espace stable de $V$ pour l'action de
$W_X$, mais également des $\tau_s, p_s$ pour $s \in \RR_X$.
L'ensemble $\RR_X$ est l'ensemble des réflexions
du sous-groupe parabolique $W_X$ de $W$, lequel est de
type $I_2(2)$ ou $I_2(3)$ car $W$ est de type ADE. Plus précisément,
tout couple $(s_1,s_2)$ de réflexions distinctes engendre un tel
sous-groupe $W_X$, de type $I_2(2)$ si $s_1 s_2 = s_2 s_1$,
de type $I_2(3)$ sinon.

\subsubsection{Si $W_X$ est de type $I_2(2)$.}
Notons $\RR_X = \{ s_1,s_2 \}$, avec $s_1 s_2 = s_2 s_1$, et
prenons $(v_{s_1},v_{s_2})$ pour base de $E_X$. Dans cette
base,
$$
\tau_{s_1} = \left( \begin{array}{cc} m & 0 \\ 0 & 1 \end{array} \right), \ \ 
\tau_{s_2} = \left( \begin{array}{cc} 1 & 0 \\ 0 & m \end{array} \right), \ \ 
$$
donc $\tau_X = (m+1) \Id$ et $[\tau_{s_1},\tau_X]=[  \tau_{s_2},\tau_X] = 0$.

\subsubsection{Si $W_X$ est de type $I_2(3)$.}
Notons $\RR_X = \{ s_1,s_2,s_3 \}$ et
prenons $(v_{s_1},v_{s_2},v_{s_3})$ pour base de $E_X$. Dans cette
base,
$$
\tau_{s_1} = \left( \begin{array}{ccc} m & -1 & -1 \\ 0 & 0 & 1 \\ 0 & 1 & 0 \end{array} \right),\ \ 
\tau_{s_2} = \left( \begin{array}{ccc} 0 & 0 & 1 \\ -1 & m & -1 \\ 1 & 0 & 0 \end{array} \right),\ \  
\tau_{s_3} = \left( \begin{array}{ccc} 0 & 1 & 0 \\ 1 & 0 & 0 \\ -1 & -1 & m \end{array} \right)
$$
donc $\tau_X = m \Id$ et $[\tau_{s_1},\tau_X]=[  \tau_{s_2} ,\tau_X]= [ \tau_{s_3},\tau_X] = 0$.
 
\subsection{Action du centre et irréductibilité}

Soit $\tau$ la somme des $\tau_s$ pour $s \in \RR$. On a
$$
\tau . v_s = \tau_s v_s + \sum_{s's \neq ss'} \tau_{s'} . v_s + \sum_{
\stackrel{s \neq s'}{s's  = ss'}} \tau_{s'}. v_s
 = m v_s + \left(\sum_{s's \neq s's} v_{s' \bullet s} - v_{s'}\right)
 + c(W) v_s$$
où l'on rappelle que $c(W) = \# \{ s' \in \mathcal{R}\setminus \{ s \}\   | \ s's = ss' \}$
ne dépend pas de $s$ puisque toutes les réflexions
sont conjuguées entre elles. De plus,
$$
\sum_{s's \neq s's} v_{s' \bullet s} - v_{s'}  = 0
$$
pour la même raison. On en déduit que $\tau$ est un scalaire,
qui vaut $m + c(W)$. Le nombre
de réflexions et $c(W)$ sont donnés par la table \ref{tablenbref}. 

\begin{table}
$$
\begin{array}{|l||c|c|c|c|c|}
\hline
W & A_{n-1} & D_n & E_6 & E_7 & E_8 \\
\hline
\# \mathcal{R} & \frac{n(n-1)}{2} & n(n-1) & 36 & 63 & 120 \\
\hline
c(W) & \frac{(n-2)(n-3)}{2} & n^2-5n+7 & 15 & 30 & 63 \\
\hline
\end{array}
$$
\caption{Dénombrement des groupes de Coxeter de type ADE}
\label{tablenbref}
\end{table}

Soit maintenant $\v \in V$ la somme des $v_s$ pour $s \in \RR$.
Il est clairement invariant par l'action de $W$, et même engendre
linéairement l'espace des vecteurs de $V$ invariants par $W$.
Pour tout $s \in \RR$,
$$
\tau_s . \v =
\tau_s. v_s + \sum_{ss' \neq s's} \tau_{s} . v_{s'} + \sum_{
\stackrel{s \neq s'}{ss'  = s's}} \tau_{s}. v_{s'}
 = m v_s + \left(\sum_{s's \neq s's} v_{s \bullet s'} - v_{s}\right)
 + \sum_{
\stackrel{s \neq s'}{ss'  = s's}}v_{s \bullet s'}$$
soit $\tau_s. \v = (m-1) v_s + s.\v - c'(W) v_s = (m-1-c'(W)) v_s + \v$
avec $c'(W) = \# \{ s' \in \mathcal{R} \ | \ s's \neq ss' \}$,
soit $c(W) + c'(W) +1 = \# \RR$.

On en déduit
\begin{proposition} \label{propirr} Si $V$ est semi-simple sous l'action des $\tau_s$,
$s \in \RR$, et
$m \not\in \{ \# \RR - c(W), -1 \}$, alors $V$ est absolument irréductible
pour cette action.
\end{proposition}
\begin{proof} Il suffit de montrer que $V$ est irréductible pour
$\k$ algébriquement clos.
Comme $m \neq -1$, l'action de $s \in \RR$ est un polynôme
de celle de $\tau_s$, donc les sous-espaces stables pour l'action
des $\tau_s$ sont stables pour l'action de $W$ puisque $\RR$ engendre
$W$. Comme l'espace des vecteurs $W$-invariants de $V$ est
de dimension 1, il existe $U \subset V$ irréductible
sous l'action des $\tau_s, s\in \RR,$ tel que $\v \in U$. Mais
alors $\tau_s . \v - \v = (m-\# \RR + c(W))v_s$ donc $v_s \in U$
pour tout $s \in \RR$, d'où $U = V$ et $V$ est irréductible.
\end{proof}

En particulier, pour $m$ générique, $V$ est irréductible
d'après la proposition \ref{fquad}. Une légère variante de
la preuve précédente permet de démontrer un résultat
d'irréductibilité de $V$ sous l'action des $\tau_s$ pour $s$
parcourant certaines parties $\SS$ de $\RR$.

\begin{proposition} \label{propirrres} Soit $\SS \subset \RR$ engendrant $W$. Si $V$ est semi-simple sous l'action des $\tau_s, s\in \SS$,
et
$m \not\in \{ \# \RR - c(W), -1 \}$, alors $V$ est absolument irréductible sous
cette action.
\end{proposition}
\begin{proof}
Comme précédemment, on suppose $\k$ algébriquement clos.
Les sous-espaces stables pour l'action
des $\tau_s, s \in \SS $ sont encore stables pour l'action de $W$ puisque $\SS$ engendre
$W$.
Comme l'espace des vecteurs $W$-invariants de $V$ est
de dimension 1, il existe $U \subset V$ irréductible
sous l'action des $\tau_s, s \in \SS,$ tel que $\v \in U$.
Soit $s \in \SS \neq \emptyset$.
On a $\tau_s . \v - \v = (m-\# \RR + c(W))v_s$ donc $v_s \in U$,
et il en est de même de tout $w. v_s = v_{wsw^{-1}}$ pour
$w \in W$ puisque $U$ est $W$-stable. Comme $W$ agit transitivement
sur $\RR$ on en déduit $U = V$ et $V$ est irréductible.
\end{proof}

\subsection{Représentation de Krammer infinitésimale}

Soit $\mathfrak{B} = \k W \ltimes \mathsf{U}\g$, où $\mathsf{U} \g$
désigne l'algèbre enveloppante universelle de $\g$. On dit qu'une
représentation $\rho$ de $\mathfrak{B}$ est \emph{in\-fi\-ni\-té\-si\-ma\-le\-ment
unitaire} par rapport à une forme quadratique si $\rho(w)$
est orthogonal pour tout $w \in W$ et $\rho(t_s)$
est autoadjoint pour tout $s \in \RR$.

\begin{proposition} \label{proprepinf}
L'application $t_s \mapsto \tau_s$ étend le $W$-module $V$ en une
représentation $\rho$ de $\mathfrak{B}$. Elle est infinitésimalement
unitaire par rapport à $(\ | \ )$ et absolument irréductible pour $m$
générique. De plus, la somme $T$ des $t_s$ pour $s \in \RR$ agit
par $m + c(W)$.
\end{proposition}

\begin{proof}
D'après ce qui précède, il reste à vérifier les relations
$[\tau_s,\tau_X]= 0$ pour tout $X$ de codimension 2 et $s \in \RR_X$.
Comme cette condition est polynomiale en $m$, il suffit de la
vérifier pour $m$ générique, et on peut donc supposer que $(\ | \ )$
est non dégénérée d'après la proposition \ref{fquad}.
Notons $s_X = \sum_{s \in \RR_X} s$. L'espace
$E_X$ est stable par $W_X$ et les $\tau_s,p_s$ pour $s \in \RR_X$,
donc il en est de même de son orthogonal $E_X^{\perp}$
par rapport à $(\ | \ )$. Comme $w p_s = p_{w \bullet s} w$
et, pour $s \in W_X$, $p_s$ est un projecteur orthogonal sur $v_s \in E_X$,
on déduit de $t_r = r - p_r$ que $[\tau_s,\tau_X] - [s,s_X]$
s'annule sur $E_X^{\perp}$, donc que $[\tau_s,\tau_X]$ s'annule
sur $E_X^{\perp}$ puisque $s s_X s = s_X$. Il suffit donc de montrer
que $[\tau_s,\tau_X]$ s'annule sur $E_X$, ce qui a déjà été fait.
L'absolue irréductibilité découle de la proposition \ref{propirr}.
\end{proof}

\subsection{Image de $\g$ dans $\gl(V)$}

Nous allons montrer que le morphisme $\g \to \gl(V)$ est surjectif.
Pour ce faire, nous établissons d'abord la proposition générale
suivante.

\begin{table}
$$
\begin{array}{|c||c|c|c|c|c|c|c|c|c|}
\hline
Type & A_r, r \geq 1 & B_r, r \geq 2 & C_r, r \geq 3 & D_r, r \geq 4 & E_6 & E_7 & E_8 & F_4
& G_2 \\
\hline
Dimension & r+1 & 2r+1 & 2r & 2r & 27 & 56 & 248 & 26 & 7 \\
\hline
\end{array}
$$\caption{Dimensions des plus petites représentations irréductibles}
\label{tableminirrep}
\end{table}

\begin{proposition} \label{propgen} Soit $V$ un $\k$-espace vectoriel de dimension
finie, et $\aa$ une sous-algèbre de Lie de $\sl(V)$ agissant
de fa\c con absolument irréductible sur $V$. Si $\rg(\aa) > \frac{1}{2} \dim
V$, alors $\aa = \sl(V)$.
\end{proposition}
\begin{proof}
On peut supposer $\k = \C$. Si $\aa$ est simple, ce résultat
découle de la classification des algèbres de Lie simples complexes
et de leur plus petite représentation irréductible (cf. table \ref{tableminirrep}).
Il reste donc à démontrer que, sous ces hypothèses, $\aa$ est simple.
Comme $\aa$ est réductive et l'algèbre dérivée $\aa'$ vérifie les mêmes hypothèses
et $\aa \subset \sl(V)$ on peut supposer que $\aa$ est semi-simple.
Décomposons $\aa = \aa_1 \times \dots \aa_p$ suivant ses facteurs simples.
On a $p \geq 1$ et on veut montrer $p =1$. Soient $r_i$ les rangs
des facteurs simples $\aa_i$. Comme $V$ est irréductible
elle s'écrit $V = V_1 \otimes \dots \otimes V_p$ avec $V_i$
représentation irréductible fidèle de $\aa_i$. En particulier,
on a $v_i = \dim V_i \geq r_i +1$. Notons $v = \dim V$, $r = \rg(\aa)$.
On a $v = v_1 \dots v_p \geq (r_1 +1)\dots (r_p+1)$ et $r > v/2$ par
hypothèse d'où $r=r_1+\dots + r_p > \frac{1}{2} (r_1 + 1)\dots ( r_p+1)$.
Si $p \geq 3$, 
$$(r_1 + 1)\dots ( r_p+1) \geqslant 1 + (r_1 + \dots + r_p)
+ r_1 r_2 + r_2 r_3 + \dots + r_{p-1}r_p + r_p r_{1} \geqslant
1 + 2r > 2r$$
ce qui contredit l'hypothèse. Il reste donc à exclure le cas $p=2$.
Notons alors $a=r_1$, $b=r_2$. On doit avoir $2(a+b) > (1+a)(1+b)$
soit $ab < a+b-1$. Or $ab \geq a+b-1$ pour tous $a,b$ entiers naturels
non nuls : supposant $a \leq b$, c'est clair si $a = 1$
et sinon $a \geq 2$ donc $ab \geq 2b \geq a+b > a+b-1$.
On en déduit $p=1$ et la conclusion.
\end{proof}

Pour $\SS \subset \RR$ on rappelle que $\g_{\SS}$ est définie comme
la sous-algèbre
de Lie de $\g$ engendrée par les $t_s$ pour $s \in \SS$.

\begin{proposition} \label{image} Pour des valeurs génériques de $m$, l'application
$\g \to \gl(V)$ est surjective. Il en est déjà de même de l'application
$\g_{\RR_0} \to \gl(V)$ où $\RR_0$ correspond aux
racines simples de $W$.
\end{proposition}
\begin{proof}
Notons
$\aa$ l'image de
$\g$ ou de $\g_{\RR_0}$ dans $\gl(V)$.
Pour $m$ générique, $V$ est semi-simple comme représentation
de $\g$ et de $\g_{\RR_0}$ d'après la proposition \ref{fquad} donc
$\aa$
est réductive et $\aa' = [\aa,\aa] \subset \sl(V)$ est
semi-simple. D'après les propositions \ref{propirr} et \ref{propirrres},
elle agit de
fa\c con irréductible sur $V$. Pour $\aa$ l'image de $\g$ il suffit de
montrer $\aa' = \sl(V)$ car $T$ agit par le scalaire $m + c(W)$, non nul
pour $m \neq - c(W)$. Comme ce scalaire est également la trace de
tout $\tau_s$ il en est de même pour l'image
de $\g_{\RR_0}$. On démontre $\aa' = \sl(V)$ par récurrence sur le rang de $W$.

Soit $\widetilde{W} \subset W$ un sous-groupe parabolique maximal
irréductible de $W$, nécessairement de type ADE,
$\widetilde{\g}$ l'algèbre de Lie et $\widetilde{\RR}$ l'ensemble
de réflexions associés, $\widetilde{\RR}_0 =\widetilde{\RR} \cap \RR_0$ et
$\widetilde{V}$ la représentation de Krammer
infinitésimale correspondante. On a un morphisme injectif
de $\widetilde{\g}$-modules $\widetilde{V} \to \Res_{\widetilde{g}} V$
induit par l'inclusion $\widetilde{\RR} \subset \RR$. En particulier,
$\aa'$ contient l'image de $[\widetilde{\g},\widetilde{\g}]$ ou
de $[\widetilde{\g}_{\widetilde{\RR}_0},\widetilde{\g}_{\widetilde{\RR}_0}]$ dans
$\sl(\widetilde{V}) \subset \sl(V)$, de rang supposé égal à $\# \widetilde{\RR} -1$.
Ainsi, $\rg(\aa') \geq \# \widetilde{\RR} - 1$ et l'on peut en déduire la
conclusion en appliquant la proposition \ref{propgen} pourvu que
l'on puisse trouver $\widetilde{W}$ avec $\#\widetilde{\RR} - 1 > \frac{1}{2}\# \RR$
et tel que $\widetilde{\RR}_0$ corresponde aux racines
simples de $\widetilde{W}$. 
Or, prenant pour parabolique maximal un groupe de type $A_{n-1},
D_{n-1},D_5,E_6,E_7$ suivant que $W$ est de type $A_n,D_n,E_6,E_7,E_8$
on constate aisément (cf. table \ref{tablenbref}) que ces hypothèses
sont vérifiées pourvu que $W$ soit de rang au moins 5. On
complète cette démonstration en vérifiant que $\g = \gl(V)$ pour
au moins une valeur de $m$, au cas par cas sur $A_1,A_2,A_3,A_4,D_4$
en prenant par exemple $m=7$ pour $A_2,A_3$ et $m=8$ pour $A_4,D_4$.
(Sur les ordinateurs actuels, le calcul effectif dure moins de dix
minutes.)
\end{proof}

\section{Identification avec les représentations de Krammer}

On suppose désormais $\k = \C$.
Nous montrons ici que la monodromie de la représentation
infinitésimale $V$ est isomorphe sur $\K = \C((h))$ à la représentation
de Krammer (généralisée) introduite dans \cite{CW,CGW}. Comme
le calcul de la monodromie est très délicat, on utilise
la rigidité de certaines algèbres quotients des algèbres de
groupe de certains groupes de réflexion (complexes) pour montrer
que ces deux représentations sont isomorphes.

On rappelle tout d'abord certaines notations classiques.
On suppose $W$ de type ADE défini par un système de racines, et on note
$\{s_1,\dots,s_n \}$ l'ensemble des réflexions simples correspondantes,
étiquetées suivant l'ordre de \cite{BOURB}. Dans le cas $ADE$ qui nous intéresse
ici, $s_i$ est relié à $s_j$ si
$<s_i,s_j>$ est de type $I_2(3)$, et ces deux réflexions ne
sont pas reliées si $<s_i,s_j>$ est de type $I_2(2)$. 

On note $R : B \to GL_N(K)$ la représentation que
l'on obtient à partir de $\rho$ par monodromie pour $\K = \C((h))$ --- on a en réalité
$R(\sigma) \in GL_N(\C[[h]])$ pour tout $\sigma \in B$. Rappelons
que l'endomorphisme $R(\sigma)$
est conjugué à $\rho(s) \exp(\ii \pi h \rho(t_s))$ si
$\sigma$ est un générateur d'Artin et $s = \pi(\sigma)$. D'autre part,
pour $m$ générique, $R$ est irréductible ainsi que sa restriction
à $P$ parce que $\rho$ est irréductible (prop. \ref{propirr}) : cela
découle des arguments de \cite{REPTHEO}, prop. 7 ou 8. Enfin, rappelons
que la classe d'isomorphisme de $R$ ne dépend pas du point base choisi.

\subsection{Factorisation par les algèbres de Hecke cubiques}

Soient $a,b,c \in  \K$. On note $H_{a,b,c}(W)$ le quotient de $\K B$
par les relations $(\ss_i -a)(\ss_i -b)(\ss_i -c)=0$ pour tout
$i$ ou encore, ce qui est équivalent en type ADE,
par la relation $(\ss_1 -a)(\ss_1 -b)(\ss_1 -c)=0$. Il est connu
que, pour $W$ de type $A_2$ et $a,b,c$ tels que les inéquations
de semi-simplicité $a \neq 0$, $a \neq b$, $a^2-ab+b^2 \neq 0$ et leurs permutés par
l'action naturelle de $\mathfrak{S}_3$ sur $a,b,c$ soient satisfaites,
ce quotient est de dimension finie et est isomorphe à l'algèbre de
groupe sur $\K$ d'un groupe de réflexions complexes de cardinal 24,
numéroté $G_4$ dans la classification de Shephard et Todd.
On renvoit à \cite{THESE,BM} pour ce résultat.

En particulier, le quotient $H_{a,b,c} = H_{a,b,c}(A_2)$ est semi-simple
et ses re\-pré\-sen\-ta\-tions ir\-ré\-duc\-ti\-bles se composent de trois représentations
de dimension 1, notées $S_x$ et définies par $\ss_i \mapsto x$ pour $x \in \{a,b,c \}$,
trois de dimension 2, notées $U_{x,y} = U_{ \{x,y \} }$ pour
$\{x,y \} \subset \{a,b,c \}$ de cardinal 2 et caractérisées par
$Sp (\ss_i) = \{x,y \}$, et une représentation de dimension 3 notée $Y$.
Dans $Y$ le spectre de $\ss_i$ est $\{a,b,c \}$. Ainsi
les représentations irréductibles de $H_{a,b,c}$ sont caractérisées
par le spectre de $\ss_1$ si $\# \{a,b,c \} = 3$. Remarquons d'autre part que
l'élément central $(\ss_1 \ss_2)^3$ agit sur $S_x$ par $x^6$, sur
$U_{x,y}$ par $-x^3y^3$ et sur $Y$ par $a^2b^2c^2$.

Notons $\rho$ la représentation de Krammer infinitésimale
associée à un groupe de Coxeter $W$ de type ADE et $R$ la représentation
de monodromie associée. Soit $\ss = \ss_1$, $s = \pi(\ss) \in W$ et
$t = t_s$. L'endomorphisme $R(\ss)$ est semi-simple
et conjugué à $\rho(s) exp( \la h \rho(t))$ avec $\la = \ii \pi$,
donc a pour spectre $q, -q^{-1}, q^m$ avec
$q = exp(\la h)$. On en déduit que cette représentation se
factorise par $H_{q, -q^{-1}, q^m}(W)$. Remarquons d'autre part que,
si $W$ est de type $A_2$, les valeurs $q, -q^{-1}, q^m$ satisfont
les inéquations de semi-simplicité pour $m \neq 1$.

\subsection{Factorisation par les algèbres de Birman-Wenzl-Murakami}

Dans \cite{BW}, Birman et Wenzl ont introduit un quotient de dimension finie
de $H_{a,b,c}(W)$ pour $W$ de type $A_n$. Cette
construction a été généralisée par Cohen, Gijsbers et Wales dans \cite{CGW}.
On peut la définir comme suit :
\begin{definition} Soit $W$ un groupe de Coxeter de type ADE.
On définit l'algèbre $BMW(W)$ sur $\Q(\alpha,l)$, où $\alpha,l$ sont des
indéterminées, comme le quotient
de l'algèbre de groupe de $B$ par les relations
$$
(\ss_1^2 + \alpha \ss_1 - 1)(\ss_1 - l^{-1}) = 0, e_i \ss_j e_i
= l e_i$$
si $s_i$ et $s_j$ sont reliées, et
où $e_i = (l/\alpha)(\ss_i^2 + \alpha \ss_i -1)$.
\end{definition}.

Pour $\alpha,l \in K\setminus \{ 0 \}$ où $K$ est un corps de caractéristique
0, on peut également définir une $K$-algèbre $BMW(W)$ spécialisée en ces
valeurs, comme quotient de l'algèbre de groupe $K B$ par les mêmes relations.
Lorsque $\alpha$ et $l$, transcendants, sont de plus algébriquement indépendants
sur $\Q$, ces deux algèbres sont isomorphes après extension des scalaires. En particulier,
prenant $K = \K$, l'algèbre $BMW(W)$ s'identifie à un quotient de
$H_{q, -q^{-1}, q^m}(W)$ pour $\alpha = q^{-1} - q$
et $l = q^{-m}$ dès que $m \not\in \Q$. En particulier, si $W = A_2$ c'est un
quotient, semi-simple, de $H_{q, -q^{-1}, q^m}$, dont les représentations
irréductibles sont $S_q$, $S_{-q^{-1}}$, $U_{q, -q^{-1}}$ et $Y$.

Soient maintenant $\ss_i, \ss_j$ tels que $s_i$ soit reliée à $s_j$ dans
le diagramme de Coxeter de $W$, et $\widetilde{B}$ le sous-groupe de
$B$ engendré par $\ss_i, \ss_j$. Soit $s_0 = s_i s_j s_i$. Le groupe
$\widetilde{B}$ est un groupe d'Artin de type $A_2$, et la restriction de $R$ à
$\widetilde{B}$ est notée $\widetilde{R}$. Pour vérifier que $R$ se factorise
par $BMW(W)$ il suffit donc de vérifier que $\widetilde{R}$ se factorise
par $BMW(A_2)$, c'est-à-dire que cette représentation
semi-simple n'admet que des composantes irréductibles de type
$S_q$, $S_{-q^{-1}}$, $U_{q, -q^{-1}}$ et $Y$, soit encore
que l'élément central $S_{ij} = (\ss_i \ss_j)^3$ n'admette pour
valeurs propres que $q^6, q^{-6}, 1$ ou  $q^{2m}$.

Notons $H_i$,
$H_j$ les hyperplans de réflexion correspondant à $\ss_i$ et $\ss_j$,
et $u_i$, $u_j$ des vecteurs
normaux à ces hyperplans. Soit
$X = H_i \cap H_j$. Décomposons $\C^n = X \oplus \C u_i \oplus \C u_j$
et notons $z_i$ resp. $z_j$ la coordonnée correspondant à
$u_i$ resp. $u_j$. On choisit $u_i$, $u_j$ de telle façon que
l'hyperplan de réflexion $H_0$ correspondant à $s_0 = s_i s_j s_i$ soit
défini par $z_i = z_j$. On note
$\alpha_H^*$ la restriction de $\alpha_H$
à $X$. Soit $\mathcal{A}' = \mathcal{A} \setminus \{ H_0,H_i,H_j \}$. On choisit $x \in X \setminus \bigcup
\mathcal{A}'$ et $\eps >0$
tel que $(x,a,b) \not\in \bigcup \mathcal{A}$ dès que $|a|,|b| \in ]0,2\eps]$.
On note $x_{\eps} = (x,\eps, 2\eps)$, qui n'appartient donc à aucun des hyperplans.
 

\begin{figure}
\includegraphics{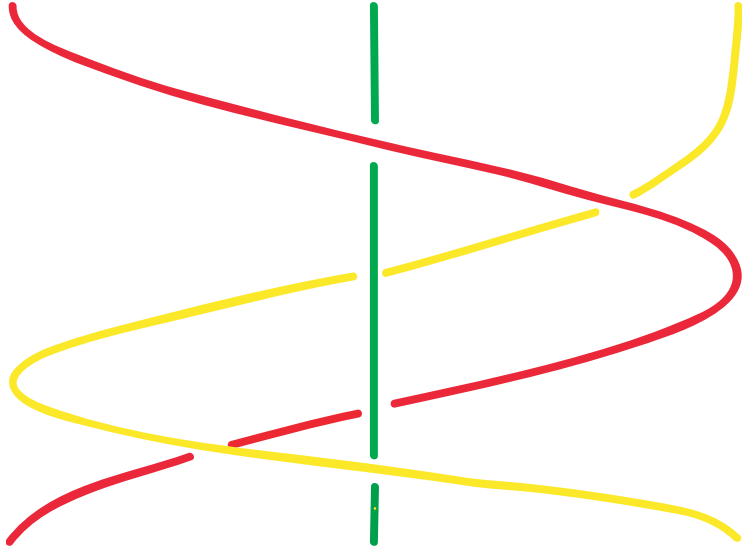}
\caption{Le lacet $S_{ij}$}
\label{tresse}
\end{figure}

Introduisons alors le lacet $x_{\eps}(t) = (x, \eps \exp(2 \ii \pi t),
2\eps \exp(2 \ii \pi t))$ de base $x_{\eps}$.
Sa classe est $S_{ij}$. Pour s'en assurer, on peut par exemple
utiliser l'identification
des triplets $(a,b,c) \in \C^3$ tels que $a,b \neq 0$, $a \neq b$ avec les
triplets $(x,y,z) \in \C^3$ tels que $x \neq y$, $y \neq z$, $x \neq z$ en
posant $a=y-x$, $b
=z-x$, $c= x+y+z$. Alors les triplets $(\eps \exp(2 \ii \pi t),
2\eps \exp(2 \ii \pi t),0)$ correspondent aux triplets
$(- \eps \exp(2 \ii \pi t), 0 , \eps \exp(2 \ii \pi t))$ qui forment
un lacet de base $(-\eps,0,\eps)$. Il est homotope au
lacet $\gamma_{\eps}(t)=(- \eps \exp(2 \ii \pi t), 0 , \eps \exp(2 \ii \pi t^2))$
par l'homotopie $(\alpha,t) \mapsto (- \eps \exp(2 \ii \pi t), 0 , \eps \exp(2 \ii \pi (\alpha t + (1-\alpha)t^2)))$.
On peut représenter ce dernier comme une tresse en le
projetant sur l'axe réel (cf. figure \ref{tresse}), ce qui permet
de l'identifier à $(\ss_i \ss_j)^3 = (\ss_j \ss_i)^3$.

Pour simplifier les notations, on note $t_k = t_{s_k}$ pour $k \in \{0,i,j\}$.
Le pullback de $\om$ sur $[0,1]$ suivant $S_{ij}$ est alors
$$
S_{ij}^* \om = h \sum_{H \in \mathcal{A}'} \rho(t_H) \frac{ O(\eps)}{\alpha_H^*(x) + O(\eps)}
+ 2 \ii \pi \rho(t_0 + t_i + t_j) dt \to 2 \ii \pi \rho(t_0 + t_i + t_j) dt
$$
quand $\eps \to 0$. Comme, lorsque $\eps$ diminue, la classe de conjugaison
de $R(S_{ij})$ ne varie pas, et que cette dernière est fermée car
$R(S_{ij})$ est semi-simple, on en déduit que $R(S_{ij})$
est conjugué à $\exp(2 \ii \pi h(t_0 + t_i + t_j))$. En particulier,
il suffit de vérifier que $\rho(t_0 + t_i + t_j)$ n'a pour valeurs propres
que $3, -3, 0$ ou $m$.

Pour ce faire, on décompose $V$ en somme de
sous-espaces stables sous l'action de $T' = t_0 + t_i + t_j$. 
Notons d'abord que $<v_{s_i},v_{s_j},v_{s_0}>$ est stable
et que $T' = t_0 + t_i + t_j$ agit sur ce sous-espace par
le scalaire $m$. On note $w_{s_k} = v_{s_k}$ pour $k \in \{0,i,j \}$.

Soit $\RR' = \RR \setminus
\{ s_i, s_j, s_0\}$. Si $u \in \RR'$
commute à chacune des trois réflexions $\{ s_i, s_j, s_0\}$,
$\C v_u$ est stable et $T'$
agit sur cette droite par la valeur $3$. On note dans ce cas $w_u = v_u$.

Si au contraire $u \in \RR'$ commute à au moins l'une d'entre elles, par exemple
$s_0$, mais pas à toutes, on montre facilement qu'elle ne
commute à aucune des deux autres. De même,
$s_i \bullet u$ ne commute alors qu'à $s_j$ et $s_j \bullet u$
ne commute qu'à $s_i$.
On note alors
$$
w_u = v_u + \frac{2}{m(m-3)} v_{s_0} + \frac{m-1}{m(m-3)} (v_{s_i} +
v_{s_j}) $$
et $w_{s_i \bullet u} = s_i \bullet w_u$,
$w_{s_j \bullet u} = s_j \bullet w_u$. On vérifie
que $t_s . w_x = w_{s \bullet x}$ pour $s \in \{ s_i,s_j,s_0\}$
et $x \in \{ u, s_i \bullet u, s_j \bullet u \}$. En particulier
$T'$ admet les valeurs propres $0$ et $3$ sur $<w_u,w_{s_i \bullet u},
w_{s_j \bullet u}> = <v_u,v_{s_i \bullet u},
v_{s_j \bullet u}>$.

On a alors défini $w_u$ pour toute
réflexion $u \in \RR'$ commutant à une seule des réflexions $s_0,s_i,s_j$.
On vérifie, facilement
en type $A_n,D_n$, et au cas par cas en types $E$,
que la dernière possibilité ($u \in \RR'$ ne commutant
à aucune des réflexions $s_i,s_j,s_0$) ne peut se produire en type ADE.
Le fait que les $w_u, u \in \RR$ forment une base de $V$
conclut : la représentation $R$ se factorise par $BMW(W)$ pour
$m$ générique, donc pour tout $m$ par analyticité de la monodromie.
\begin{proposition} \label{propidentif} Pour tout $m \in \C$, la représentation $R$ se factorise
par $BMW(W)$ spécialisée en $\alpha = q^{-1}-q$, $l = q^{-m}$.
\end{proposition}

\subsection{Annulation sur l'idéal $I_2$}

Dans \cite{CGW} est introduit un idéal remarquable de l'algè\-bre $BMW(W)$,
noté $I_2$, qui est engendré par les éléments
$e_i e_j$ pour $s_i$ non reliée à $s_j$, c'est-à-dire si $s_i$ et
$s_j$ commutent. Nous allons montrer que $R$ s'annule sur
$I_2$, c'est-à-dire $R(e_i e_j) = 0$. Cette fois encore, on peut supposer
$m$ générique.

Soient $s_i$ et $s_j$ non reliées. On note $H_i$, $H_j$ les
hyperplans de réflexion correspondants, et on peut décomposer $\C^n =
X \oplus \C u_i \oplus \C u_j$ avec $X = H_i \cap H_j$,
$s_i$ resp. $s_j$ agissant par $-1$ sur $u_i$ resp. $u_j$. On note
$z_k$ la coordonnée correspondant à $u_k$ pour $k \in \{ i,j \}$,
et $\alpha_H^*$ la restriction de $\alpha_H$ à $X$.
Il
existe $x \in X$ tel que $\alpha_H^*(x) \neq 0$ pour
tout $H \in \mathcal{A}' = \mathcal{A} \setminus \{ H_i ,H_j \}$. On note
$x_{\eps} = (x, \eps, \eps)$. Il existe $\eps_0 >0$ tel que,
pour tout $\eps \in ]0,\eps_0[$ et pour tout $H \in \RR$ on
ait $\alpha_H(x_{\eps}) \neq 0$, et même $\alpha_H((x,u \eps, v \eps)) \neq 0$
pour tous $u,v$ de module 1. Prenant un tel $x_{\eps}$ pour point
base, un chemin représentant $\ss_i$ (resp. $\ss_j$)
est $\ss_{i,\eps} = (x, \eps \exp(2 \ii \pi t), \eps)$
(resp. $\ss_{j,\eps} = (x, \eps,\eps \exp(2 \ii \pi t))$). On
note $\tau_i = \tau_{H_i}$ avec $\tau_H = \rho(t_H)$. On a alors
$$
\frac{1}{\ii \pi} \ss_i^* \om = h \left( \tau_{i} + \sum_{H
\in \mathcal{A}'} \frac{O(\eps)}{\alpha_H^*
(x) + O(\eps)} \tau_H \right) dt,
\frac{1}{\ii \pi} \ss_j^* \om = h\left( \tau_{j} + \sum_{H
\in \mathcal{A}'} \frac{O(\eps)}{\alpha_H^*
(x) + O(\eps)} \tau_H \right) dt
$$
d'où $\ss_i^* \om \to h \ii \pi \tau_{i} dt$
et $\ss_j^* \om \to h \ii \pi \tau_{j} dt$ quand $\eps \to 0$.

En particulier,
si $R_{\eps}$ désigne la monodromie de $\rho$ avec point base $x_{\eps}$,
$R_{\eps}(\ss_i) \to \sss_i = \rho(s_i) \exp(\ii \pi h \tau_i)$
et
$R_{\eps}(\ss_j) \to \sss_i = \rho(s_j) \exp(\ii \pi h \tau_j)$
quand $\eps \to 0$.
Notons maintenant $E = e_i e_j$. On a $e_k^2 = \alpha e_k$
pour $k \in \{i,j \}$ et $e_i e_j = e_j e_i$,
donc $E^2 = \alpha^2 E$. Ainsi,
$R(E) = 0$ ssi $tr( R(E)) = 0$. Plus précisément, le spectre
de $E$ est inclus dans $\{ 0, \alpha^2 \}$. Comme
$\alpha^2 \equiv 
-4 \pi^2 h^2 \mod h^3$ on a donc $R(E) = 0$ ssi $tr(R(E)) \equiv 0 \mod h^3$.

Remarquons
que, puisque $R$ se factorise par $BMW(W)$, on a nécessairement
$tr(R(E)) \equiv 0 \mod h^2$. D'autre part, les endomorphismes
$R_{\eps}(E)$ sont conjugués pour $\eps \in ]0,\eps_0[$ donc
$tr(R_{\eps}(E))$ ne dépend pas de $\eps$. Il suffit alors de montrer que
$tr(R_{\eps}(E))$ modulo $h^3$ tend vers 0 quand $\eps \to 0$,
donc que $\ee_i \ee_j \equiv 0 \mod h^3$ avec
$\ee_k = (l/\alpha)( \sss_k^2 + \alpha \sss_k -1)$
pour $k \in \{i,j \}$. Or on a
$$
\ee_k =\rho \left( (s_k - t_k)(1- (m-t_k) \tilde{h}) + \tilde{h}^2 X_k + O(\tilde{h}^3) \right)
$$
avec
$$
X_k = \frac{m^2}{2}(s_k - t_k) +m t_k(t_k - s_k) + \frac{t_k}{6}
(-4 t_k^2 + 3 s_k t_k +1) , \ \ \tilde{h} =  \ii \pi h
$$

Remarquons ensuite que, puisque $s_i$ et $s_j$ ne sont pas
reliées, $s_i - t_i$ et $s_j - t_j$ commutent et leurs images par
$\rho$ sont, à un scalaire
près et pour des valeurs génériques de $m$, des projecteurs orthogonaux sur
des droites orthogonales pour $(\ | \ )$. On en déduit que
le produit de ces images est nul, ceci pour toute valeur de $m$.
Ainsi,
$$
\ee_i \ee_j \equiv \tilde{h}^2\rho\left(  (s_i - t_i) X_j + X_i(s_j - t_j) \right) 
\mod h^3.
$$
et, notant $x \sim y$ si $\rho(x) = \rho(y)$,  
$(s_i - t_i)X_j \sim  
(s_i - t_i)\frac{t_j}{6}
(-4 t_j^2 + 3 s_j t_j +1)$. De plus,
$(s_i - t_i)t_j \sim (s_i - t_i)s_j$,
d'où
$(s_i - t_i)X_j \sim  
(s_i - t_i)\frac{s_j}{6}
(-4 s_j^2 + 3 s_j s_j +1) = (s_i - t_i) \frac{s_j}{6}(-4 + 3 + 1) = 0 $
et $\ee_i \ee_j \equiv 0 \mod h^3$.
On a ainsi montré

\begin{proposition} Pour tout $m \in \C$, la représentation $R$ se factorise par $BMW(W)/I_2$
pour $\alpha = q^{-1}-q$, $l = q^{-m}$.
\end{proposition}

\subsection{Identification à la représentation de Krammer}

Dans \cite{CGW} la structure de $BMW(W)/I_2$ est analysée lorsque $\alpha$ et $l$ sont
algébriquement indépendants. On suppose donc $m \not\in \Q$. Il existe
un sous-groupe parabolique $\widetilde{W}$ de $W$ telle que l'algèbre de Hecke
$\widetilde{H}$ de $\widetilde{W}$ avec paramètres $q$, $-q^{-1}$ détermine
toutes les représentations irréductibles de $BMW(W)/I_2$. Plus précisément,
à toute (classe de) représentation irréductible $\theta$ de $\widetilde{H}$
est associée une représentation irréductible $\Gamma_{\theta}$ de
$BMW(W)/I_2$ de dimension $(\# \RR) \dim(\theta)$, et
cette correspondance $\theta \to \Gamma_{\theta}$ entre classes de représentations
irréductibles des deux algèbres est bijective.

A tous $u \in \RR$ et $i \in \{1,\dots,n\}$ tel que $u \neq s_i$ et
$us_i = s_i u$ les auteurs de \cite{CGW} associent un élément $h_{u,i} \in \widetilde{H}$ tel
que la trace de $\Gamma_{\theta}(\sigma_i)$ soit
$$
q^m + \sum_{\stackrel{u \neq s_i}{us_i = s_i u}} \theta(h_{u,i}) 
+ \frac{q-q^{-1}}{2} c'(W).
$$
La représentation de Krammer généralisée
est alors $\Gamma_{\theta}$ avec $\theta$ entièrement déterminé par
$\theta(h_{u,i}) = q$. D'autre part, on a toujours $(h_{u,i}-q)(h_{u,i}+
q^{-1}) = 0$, donc elle est caractérisée parmi les représentations
$\Gamma_{\theta}$ par le fait d'être de dimension $\# \RR$ et d'envoyer
$\sigma_1$ sur un élément de trace $q^m + q c(W) + (q-q^{-1})\frac{c'(W)}{2}$.
En effet, si une telle représentation est de dimension $\# \RR$
on a $\dim \theta = 1$ donc $\theta(h_{u,i}) \in \{q,q^{-1} \}$, et
$$
\sum_{\stackrel{u \neq s_i}{us_i = s_i u}} \theta(h_{u,i}) =
q c(W) \Leftrightarrow \forall u,i\ \ \theta(h_{u,i}) = q.
$$

Comme, pour $m$ générique, $R$ est une représentation irréductible de $BMW(W)/I_2$ de dimension
$\# \RR$, il suffit donc de montrer que $tr(R(\sigma_1))$ a la valeur
voulue. Comme $R(\sigma_1)$ est conjugué à $\sss = \rho(s \exp(\la t_s))$,
il suffit de calculer la trace de $\sss$.

On a $V = \C v_s \oplus U \oplus U'$ avec $U = \bigoplus_{\stackrel{u \neq s}{us = s u}}
 \C v_u$ et $U' = \bigoplus_{us \neq su} \C v_u$ ; de plus $\dim(U) = c(W)$,
$\sss. v_s = q^m v_s$, $\sss . v_u  = q v_u$ si $u \neq s$ et $us = su$.
D'autre part $U'$ est somme directe de $\frac{c'(W)}{2}$ plans
stables pour l'action de $s$ et $t_s$, engendrés par les familles
de la forme $(v_u,v_{u'})$ telles que $<s,u,u'>$ soit de type $I_2(3)$. Sur
ces
plans, $s= t_s$ admet pour vecteurs propres $v_u + v_{u'} +
\frac{2}{m-1} v_s$, $v_u - v_{u'}$ associés respectivement aux valeurs
propres $1$ et $-1$. Ainsi, la restriction de $\sss$ à un de ces plans
a pour trace $q-q^{-1}$, et on en déduit que $tr(R(\sigma_1))$ a la
valeur voulue d'où le

\begin{theoA} 
Pour tout $W$ de type $ADE$
et $m \not\in \Q$,
la représentation de Krammer $R_{\mathrm{K}}$ plongée dans $GL_N(\K)$
par $q = e^h$ et $q= e^{m h}$ est isomorphe à $R$.
\end{theoA}

\begin{ccor} Les groupes d'Artin purs de type sphérique sont résiduellement
nilpotent-sans-torsion. En particulier, ils sont bi-ordonnables.
\end{ccor}
\begin{proof} On peut supposer que $W$ est un groupe
de Coxeter irréductible. D'après \cite{CRISP}, tout
groupe d'Artin pur de type sphérique se plonge
dans un groupe d'Artin pur de type ADE. On peut donc
supposer que $W$ est de type ADE. Notons toujours $\mathcal{M} = \C[[h]]$,
$(x,y) = xyx^{-1}y^{-1}$ et $C^1P = (P,P)$, $C^{r+1} P = (P,C^rP)$
les termes de la suite centrale descendante de $P$. Il s'agit de
montrer en particulier que leur intersection est triviale. Identifiant $V$
à $\C^N$, on a $R(P)
\subset 1 + h M_N(\mathcal{M})$ donc $R(C^r P) \subset 1 + h^{r+1} M_N(\mathcal{M})$
d'où
$$R\left(\bigcap_{r\geq 1} C^r P \right) = \{1 \} \Rightarrow
\bigcap_{r\geq 1} C^r P = \{ 1 \}$$
par fidélité de $R$, au moins si les paramètres $q$ et $t$ de \cite{CW}
sont transcendants sur $\C$ et algébriquement indépendants, ce qui est
le cas dès que $m \not\in \Q$. Autrement dit, $R$ plonge $P$ dans le groupe $1+h M_N(\mathcal{M})$ qui
est non seulement résiduellement nilpotent mais aussi résiduellement nilpotent-sans-torsion,
donc bi-ordonnable (cf. par exemple \cite{PASS}).
\end{proof}

Remarquons que les modèles matriciels de la représentation
de Krammer per\-met\-tent éga\-le\-ment d'obtenir ce dernier résultat,
en posant un peu artificiellement $q = e^h$ et $t = e^{m h}$
pour $m \not\in \Q$ dans les modèles
de \cite{CW} et en réduisant les formules modulo $h$. Nous renvoyons
à \cite{resnil} où nous avons utilisé cette méthode.

\section{Enveloppe algébrique et unitarité de la représentation de Krammer}

\subsection{Enveloppe algébrique}

Le calcul de l'image de $\rho$ nous permet de déterminer l'enveloppe
algébrique des groupes d'Artin de type ADE dans la représentation
de Krammer associée. Notons $V_h = V \otimes \K$.
Pour toute réflexion $s$ de $W$, il existe
un $\sigma \in P$ tel que $R(\sigma)$ est congru à $exp(2 \ii \pi h \rho(t_s))$
modulo $h^2$ ($\sigma$ est par exemple un lacet faisant un tour
complet autour de l'hyperplan associé à $s$).

On déduit alors
des arguments de \cite{REPTHEO}, prop. 24
(ou bien \cite{LIETRANSP} section 7.1) que l'algèbre de Lie
de l'enveloppe algébrique de $R(B)$ contient l'image de $\g$
tensorisée par $\K$. D'après la proposition \ref{image}, on sait
que c'est $\gl_N(\k) \otimes \K = \gl_N(\K)$. On en déduit 

\begin{theoB} 
Pour tout groupe de Coxeter fini irréductible de type ADE et $m \not\in \Q$,
la représentation $R$ identifie $B$ à un sous-groupe Zariski-dense
de $GL_N(\K)$, où $N$ est le nombre de réflexions de $W$.
Il en est de même des images des sous-groupes $C_W^k$
engendrés par les éléments $\sigma_1^{2k},\dots, \sigma_n^{2k}$
où les $\sigma_i$ désignent les générateurs d'Artin de $B$
et $k \geq 1$.
\end{theoB}

Notons $\LL$ le corps des fractions de $\mathcal{L}$.
Comme, pour tout $m \not\in \Q$,
la représentation $R$ est isomorphe à la représentation de Krammer
$R_{\mathrm{K}}$ après extension des scalaires,
on déduit immédiatement de ce théorème le

\begin{ccor}
Pour tout groupe d'Artin sphérique irréductible de type ADE, la
représentation de Krammer $R_{\mathrm{K}}$
identifie $B$ à un sous-groupe Zariski-dense de $GL_N(\LL)$,
\end{ccor}

\subsection{Remarques sur d'autres sous-groupes}

Ce théorème de Zariski-densité implique, par des arguments généraux que
nous détaillons plus loin, la Zariski-densité de nombreux
sous-groupes de ces groupes : essentiellement des sous-groupes
sous-normaux et les sous-groupes d'indice fini.
En plus de ceux-ci, on peut également raffiner la démonstration de ce théorème
afin d'obtenir la Zariski-densité d'autres sous-groupes qui ne
rentrent pas dans ce cadre. Nous donnons ici un exemple en
type $A_{n-1}$. Rappelons que, en type $A_{n-1}$, le groupe d'Artin
pur est engendré par les éléments 
$$\xi_{ij} = \ss_{j-1}\dots \ss_{i+1}
\ss_i^2 \ss_{i+1}^{-1} \dots \ss_{j-1}^{-1}$$
pour $1 \leq i<j \leq n$, et les sous-groupes $F_n^k = < \xi_{in}^k \ | \
1 \leq i < n >$ pour $k \geq 1$ sont des sous-groupes remarquables
de $B$, qui ne sont sous-normaux que si $k = 1$.

\begin{proposition}
Soit $W$ de type $A_{n-1}$, et $m \not\in \Q$. Alors,
pour tout $k \geq 1$, le groupe $R(F_n^k)$ est Zariski-dense dans $GL_N(\K)$.
\end{proposition}

La démonstration est identique, sous la réserve de montrer que
que l'application
$\g_{\SS} \to \gl(V)$
pour $\SS = \{(1 \ n), \dots, (n-1\ n) \}$ est surjective.
La démonstration de ce dernier fait est similaire à celle de la proposition \ref{image}.
L'hypothèse non vérifiée est que $\SS \cap \mathfrak{S}_{n-1}$ est vide.
Pour résoudre
ce problème, on remarque que l'image de $\g_{\SS}$ est conjuguée
à l'image de $\g_{\SS'}$ avec $\SS' = \{ (1\ n-1), \dots, (n-2\ n-1), (n-1 \ n) \}$ par
l'image de la transposition $(n\!-\!1\ n)$. On montre alors que l'image de
$\g_{\SS'}$
est $\gl(V)$ par le même argument puisque $\SS' \cap \mathfrak{S}_{n-1}$ engendre
$\mathfrak{S}_{n-1}$ et il en est donc de même de
l'image de $\g_{\SS}$.

\subsection{Unitarité}

Nous discutons ici du lien entre l'unitarité infinitésimale
de $\rho$ et l'éventuelle unitarisabilité de $R$. Soit $\k$
un corps de caractéristique 0 et supposons $\g$, $\mathfrak{B}$
et $\rho$ définies sur $\k$. Les algèbres $\g$ et
$\mathfrak{B}$ sont naturellement graduées, par $\deg t_s = 1$
pour tout $s \in \RR$, et $\deg(w) = 0$ si $w \in W$.
On note $\widehat{\mathfrak{B}}$ et $\widehat{\g}$ les complétions de
$\mathfrak{B}$ et $\g$
par rapport à cette graduation.
On a $\exp \widehat{\g} \subset
\widehat{\mathfrak{B}}$. Pour déduire de $\rho$ une
représentation de $B$ sur $\K = \k((h))$ il suffit
d'avoir un morphisme $\PPhi : B \to \widehat{\mathfrak{B}}$.
En effet, notant $\tilde{\rho}$ la représentation de
$\widehat{\mathfrak{B}}$ dans $\End(V_h)$ définie par
$\tilde{\rho}(w) = \rho(w)$ pour $w \in W$ et $\tilde{\rho}(t_s)
= h \rho(t_s)$, $\widehat{\Phi}(\rho) = \tilde{\rho} \circ
\PPhi$ définit une représentation de $B$ dans $GL(V_h)$.
Un tel morphisme est donné pour $\k = \C$ par monodromie,
et permet d'obtenir la représentation de Krammer $R$.

Supposons que l'on ait un morphisme $\PPhi : B \to \widehat{\mathfrak{B}}$
satisfaisant la condition suivante
\begin{eqnarray}
\exists \Phi_1,\dots,\Phi_n \in \exp \widehat{\g} \ | \ \PPhi(\ss_i)
= \Phi_i s_i \exp( t_{s_i}) \Phi_i^{-1}
\label{condfond}
\end{eqnarray}

Alors la représentation $\widehat{\Phi}(\rho)$
correspondante de $B$ prend ses valeurs dans le groupe
unitaire formel $U_N^{\varepsilon}(K)$ correspondant à l'automorphisme
involutif $\varepsilon : f(h) \mapsto f(-h)$ de $K$ et à la forme
$(\ | \ )$. Si $\k \subset \R$, on a montré dans \cite{REPTHEO} que,
après torsion des coefficients de la représentation par un
automorphisme convenable de $K$ de manière à faire converger les
coefficients pour $h \in \ii \R$ suffisamment proche de 0,
on pouvait en déduire des représentations unitaires $\tilde{R}$ de $B$, pourvu
que $(\ | \ )$ soit définie positive. On a montré ici que tel était le
cas pour $m >  \# \RR - c(W)$.

L'existence de tels morphismes satisfaisant (\ref{condfond}) est établie
si $B$ est \og fiber-type \fg\ au sens de Falk et Randell, c'est-à-dire
si $W$ est de type $A_n$, $B_n$ ou diédral. C'est en effet une
conséquence des travaux de Drinfeld \cite{DRINFELD} en type $A_n$, Enriquez
\cite{ENRIQUEZ} en
type $B_n$ et de \cite{DIEDRAUX} pour les groupes diédraux.
Malheureusement ce résultat reste pour l'instant à l'état de conjecture
en types $D_n$ et $E_6,E_7,E_8$.

En type $A_n$, cela nous permet en tous cas de retrouver abstraitement
un résultat de R. Budney, qu'il démontre dans \cite{BUDNEY} en construisant
une forme quadratique explicite.
\begin{proposition} \label{propunita} La représentation de Krammer en type $A_{n-1}$
est unitarisable pour $|q| = 1$ proche de 1 et $m > \# \RR - c(W) = 2n-3$.
\end{proposition}
\begin{proof}
Nous utilisons les notations de \cite{REPTHEO}. Soit $\Phi \in
\mathbbm{Assoc}_1(\Q)$ un associateur rationnel et $\PPhi$ le morphisme
de Drinfeld associé. On démontre
facilement que $\widehat{\Phi}(\rho)$ se factorise par $BMW(A_n)$ :
cette représentation irréductible se factorise par l'algèbre de
Birman-Wenzl-Murakami soit pour les mêmes raisons que la représentation $R$,
soit en utilisant les arguments de \cite{QUOT} prop. 4 ; l'identification
avec la représentation de Krammer découle alors
du fait que $\widehat{\Phi}(\rho)(\delta_r) = \exp(2h \rho(Y_r))$,
avec $Y_r = t_{1r}+\dots+t_{r-1,r}$, $\delta_r = \sigma_{r-1} \dots \sigma_2
\sigma_1^2 \sigma_2 \dots \sigma_{r-1}$, et que le spectre de
l'action des $\delta_r$ permet d'identifier la représentation
irréductible de l'algèbre de Birman-Wenzl-Murakami considérée (cf. \cite{QUOT},
section 9). Après torsion de la représentation par
un automorphisme adéquat du corps des coefficients
et spécialisation en $h \in \ii \R$ suffisamment proche de 0
on en déduit des représentations unitaires $\tilde{R}_h$ (cf. \cite{REPTHEO},
appendice A). Pour des valeurs
génériques de $h$, ces représentations sont bien isomorphes
aux représentations de Krammer au moins si l'on a pris soin de choisir
l'automorphisme de telle façon que l'image des générateurs d'Artin soit suffisamment
proche, modulo $h$, de leur image dans la représentation originelle.
\end{proof}

\section{Applications aux groupes d'Artin-Tits}

\subsection{Zariski-densité : preuve du théorème 4}

Le but de cette section est de démontrer le théorème 4, à savoir que
tout groupe d'Artin-Tits sphérique et irréductible se plonge
de façon Zariski-dense dans un groupe linéaire $GL_N(K)$.
On peut prendre pour $\k$ un corps de caractéristique 0
quelconque et $K = \K$. Pour simplifier l'exposé, on
supposera $\k = \C$.
D'après le théorème B, c'est le cas lorsque $W$
est de type $ADE$. Lorsque $B$ est de type $B_n$, il se
plonge dans un groupe de type $A_{n}$. En effet, si
l'on note $\ss_1,\dots,\ss_{n-1}$ les générateurs d'Artin naturels
d'un groupe $G$ de type $A_{n-1}$, il est classique que le sous-groupe
$H$ engendré par $\ss_1^2,\ss_2,\dots,\ss_n$ est de type $B_n$, et
est d'indice 2 dans $G$. De plus, il contient le sous-groupe
$C_n^1$ de $G$. On en déduit que, notant $R : G \to GL_N(\K)$
la représentation construite en type $A_{n-1}$, $R(H)$
est également Zariski-dense dans $GL_N(K)$, donc le résultat est vrai
en type $B_n$.

Pour traiter les groupes de type $I_2(m)$ nous allons utiliser la représentation
de Burau. Rappelons que, pour tout groupe de Coxeter $W$, la représentation de Burau réduite est une déformation
de la représentation de réflexion de $W$. Plus précisément,
il s'agit de la représentation (irréductible) de l'algèbre de Hecke de $W$ qui correspond à
la représentation de réflexion de $W$. On définit ici l'algèbre de Hecke de $W$ comme
le quotient de l'algèbre de groupe $\K B$ par l'idéal engendré par
les éléments $(\ss - q)(\ss + q^{-1})$, où $q = e^h$ et $\ss$
parcourt les générateurs d'Artin de $B$.

Si $\rho_0 : W \to GL_n(\C)$ est la
représentation de réflexion de $W$, alors il existe
une re\-pré\-sen\-ta\-tion équivariante $\rho : \g \to \gl_n(\C)$ définie par $\rho(t_s) = \rho_0(s)$
pour tout $s \in \RR$, et il est classique que la monodromie $R_{\mathrm{bur}} : B \to GL_N(\K)$ d'une telle représentaton
est isomorphe à la représentation de Burau réduite de $B$. Définissont $\tilde{R}_{\mathrm{bur}} : B \to GL_N(\K)$
par $\tilde{R}_{\mathrm{bur}}(\ss) = q R_{\mathrm{bur}}(\ss)$ si $\ss$ est un générateur d'Artin. Alors $\tilde{R}_{\mathrm{bur}}$
est la monodromie de $\tilde{\rho} : t_s \mapsto \Id + \rho_0(s)$.

Revenant au cas des groupes de type $I_2(m)$ pour $m \geq 2$, il a été démontré par G. Lehrer et N. Xi que
$\tilde{R}_{\mathrm{bur}}$ était fidèle (cf. \cite{LEHRER}, prop. 4.1). Pour conclure comme précédemment
à la Zariski-densité de $\tilde{R}_{\mathrm{bur}}(B)$ et de ses sous-groupes $\tilde{R}_{\mathrm{bur}}(C_W^k)$,
il suffit de vérifier que, si une présentation de $W$ est donnée par
$<s,t \ | \ s^2=t^2=1, (st)^m=1 >$, alors les matrices
$$
\tilde{\rho}(s) = \left( \begin{array}{cc} 0 & 0 \\ -c & 2 \end{array} \right)
\ \ \ \ \ \ \tilde{\rho}(t) = \left( \begin{array}{cc} 2 & -c \\ 0 & 0 \end{array} \right)
$$
avec $c = 2\cos(\frac{2\pi}{m})$, engendrent l'algèbre de Lie $\gl_2(\C)$,
ce qui se vérifie immédiatement.

Il reste à démontrer ce résultat pour les groupes exceptionnels de type $F_4,H_3$ et $H_4$. Nous numérotons
les générateurs d'Artin $\ss_i$ de ces groupes exceptionnels comme suit :
$$
\begin{array}{|c|c|c|}
\hline
F_4 & H_3 & H_4 \\
\xymatrix{
1 \ar@{-}[r] & 2 \ar@2{-}[r] & 3 \ar@{-}[r] & 4
} 
& \xymatrix{
1 \ar@3{-}[r] & 2 \ar@{-}[r] & 3  
} 
& \xymatrix{
1 \ar@3{-}[r] & 2 \ar@{-}[r] & 3 \ar@{-}[r] & 4 
} \\
\hline
\end{array}
$$

Ces trois groupes admettent des plongements naturels dans des groupes de type $ADE$, plus précisément
dans les groupes suivants
$$
\begin{array}{|c|c|}
\hline
E_6 & D_6 \\
\xymatrix{
1 \ar@{-}[r] & 3 \ar@{-}[r] & 4 \ar@{-}[r] \ar@{-}[d] & 5 \ar@{-}[r] & 6  \\
 & & 2  
} &
\xymatrix{
1 \ar@{-}[dr] \\
 &  3 \ar@{-}[r] & 4 \ar@{-}[r] & 5 \ar@{-}[r] & 6 \\
2 \ar@{-}[ur]
}  \\
\hline
\multicolumn{2}{|c|}{ E_8 } \\
\multicolumn{2}{|c|}{ \xymatrix{
1 \ar@{-}[r] & 3 \ar@{-}[r] & 4 \ar@{-}[r] \ar@{-}[d] & 5 \ar@{-}[r] & 6 \ar@{-}[r] & 7 \ar@{-}[r] & 8 \\
 & & 2  
}} \\
\hline
\end{array}
$$
et les plongements des premiers dans les seconds sont
définis comme suit
$$
\begin{array}{|lcl|rcl|rcl|}
\hline
\mbox{\og} F_4 & \to & E_6 \mbox{\fg} & \mbox{\og}H_3 & \to & D_6 \mbox{\fg}& \mbox{\og}H_4 & \to & E_8 \mbox{\fg}\\
\hline
\ss_1 & \mapsto & \ss_2 & \ss_1 & \mapsto & \ss_4 \ss_1 & \ss_1 & \mapsto & \ss_2 \ss_5 \\ 
\ss_2 & \mapsto & \ss_4 & \ss_2 & \mapsto & \ss_5 \ss_3 & \ss_2 & \mapsto & \ss_4 \ss_6 \\ 
\ss_3 & \mapsto & \ss_3 \ss_5 & \ss_3 & \mapsto & \ss_6 \ss_2 & \ss_3 & \mapsto & \ss_3 \ss_7 \\ 
\ss_4 & \mapsto & \ss_1 \ss_6 &  &  &  & \ss_4 & \mapsto & \ss_1 \ss_8 \\ 
\hline
\end{array} 
$$
Ces plongements font partie des \og homomorphismes de pliage \fg\ définis
par J. Crisp dans \cite{CRISP}. Ils sont injectifs, et préservent l'élément de Garside $w_0$.
Cet élément de $B$ engendre le centre en types $F_4,H_3,H_4,D_6,E_7,E_8$.
En type $E_6$ en revanche, le centre est engendré par $w_0^2$. Cela implique
que la représentation de Krammer restreinte $R_{\mathrm{K}} \circ \varphi$ n'est
pas absolument irréductible en type $F_4$, et en particulier ne peut être
d'image Zariski-dense. Plus précisément, l'action de l'élément de Garside
du groupe d'Artin de type $F_4$ décompose la représentation de Krammer de type $E_6$, de dimension 36,
en somme de deux composantes, respectivement de dimension 24 et 12 (on
note que le groupe de Coxeter de type $F_4$ admet 24 réflexions). F. Digne a montré dans \cite{DIGNE} que l'action de $F_4$
sur cette composante
de dimension 24 était fidèle. Cela signifie que la représentation correspondante $R_{\mathrm{K}'} : B \to GL_{24}(K)$
pour $B$ de type $F_4$ est injectif. Pour $B$ de type $H_3$ ou $H_4$, on notera
$R_{\mathrm{K}'} = R_{\mathrm{K}} \circ \varphi$.

Il reste à vérifier que les $\log R_{\mathrm{K}'}(\ss^2)$ engendrent
$\gl_N(K)$, pour $\ss$ parcourant les générateurs d'Artin de $B$, quand
$B$ est de type $F_4,H_3,H_4$, et $N = 24,30,120$, respectivement. Pour
ce faire, il suffit de vérifier que les $ u(\ss) = (h^{-1}\log R_{\mathrm{K}'}(\ss^2)) \mod h$
engendrent $\gl_N(\k)$, c'est-à-dire qu'ils engendrent une algèbre de Lie de dimension
$N^2$, au moins pour une certaine valeur de $m$. On obtient des modèles
matriciels\footnote{\samepage
Dans \cite{DIGNE} il convient de corriger les deux données suivantes. Page 56,
si l'on note $\alpha_i$ la racine simple correspondant à $s_i$,
notamment $\alpha_2<\alpha_3$, l'ordre choisi pour les racines
positives $x \alpha_1 + y \alpha_2 + z \alpha_3 + t \alpha_4$ n'est
pas celui indiqué. La liste ordonnée des quadruplets $[x,y,z,t]$ est au contraire
$$
\begin{array}{c}
 {[} 1, 0, 0, 0 {]}, {[} 0, 0, 0, 1 {]}, {[} 0, 1, 0, 0 {]}, {[} 0, 0, 1, 0 {]}, {[} 1, 1, 0, 0 {]}, {[} 0, 0, 1, 1 {]}, {[} 0, 1, 1, 0 {]}, 
 {[} 1, 1, 1, 0 {]}, \\  {[} 0, 1, 1, 1 {]}, {[} 0, 2, 1, 0 {]}, {[} 1, 1, 1, 1 {]}, {[} 1, 2, 1, 0 {]}, {[} 0, 2, 1, 1 {]}, {[} 1, 2, 1, 1 {]}, 
  {[} 2, 2, 1, 0 {]}, {[} 0, 2, 2, 1 {]}, \\
{[} 1, 2, 2, 1 {]}, {[} 2, 2, 1, 1 {]}, {[} 1, 3, 2, 1 {]}, {[} 2, 2, 2, 1 {]}, {[} 2, 3, 2, 1 {]}, 
  {[} 2, 4, 2, 1 {]}, {[} 2, 4, 3, 1 {]}, {[} 2, 4, 3, 2 {]} 
\end{array}
$$
D'autre part, les formules décrivant la représentation en type $D_n$
sont incomplètes. Dans la description de l'image de $\rho(s)(e_r)$,
page 52, il faut corriger le quatrième cas en
$$
\begin{array}{ll}
q e_{srs} + (1-q) e_r + t\frac{(q-1)^2(q^i-1)}{q^i} & \mbox{ si }
\alpha+\beta \in \Phi, n(\alpha+\beta) = i \geq 1, n(\beta) = i-1 \\
q e_{srs} + (1-q) e_r  & \mbox{ si }
\alpha+\beta \in \Phi, n(\alpha+\beta) = i \geq 1, n(\beta) \neq i-1 \\
\end{array}
$$
} à coefficients entiers à partir de \cite{DIGNE}, donc on peut
supposer $u(\ss) \in \gl_N(\Z)$. Le calcul 
sur ces matrices à coefficients entiers a été fait informatiquement,
à l'aide d'un programme en langage C et en réduisant modulo 19 pour augmenter la vitesse
d'exécution et réduire la taille mémoire nécessaire. En choisissant par exemple les valeurs $m=5$ pour les types $H_3$ et $F_4$,
$m=7$ pour les types $H_3$ et $H_4$, on obtient que l'algèbre de Lie engendrée
par les matrices $\overline{u(\ss)} \in \gl_{N}(\F_{19})$ est $\gl_N(\F_{19})$,
donc l'algèbre de Lie engendrée par les $u(\ss)$ est de dimension $N^2$
sur $\Q$ donc sur $\k$ pour au moins une valeur de $m$. Ceci conclut la
démonstration du théorème 4.

\subsection{Application aux décompositions en produit direct}

Nous appliquons ici le résultat de Zariski-densité de la représentation
de Krammer (généralisée) au problème de la dé\-com\-po\-si\-tion en produit
direct des groupes d'Artin et de ses sous-groupes remarquables.
Ces résultats peuvent également se déduire de la section suivante.
Nous en donnons ici une preuve directe de nature différente.
Le lemme clé est le suivant.

\begin{lemme} Soit $G$ un groupe, et $V$ une représentation
semi-simple fidèle de $G$ sur un corps $K$ algébriquement clos de caractéristique 0
telle que $V \otimes V$ ait au plus trois composantes
irréductibles. Si $G$ est produit direct $G_1 \times G_2$
de deux de ses sous-groupes, alors $G_1 \subset Z(G)$ ou
$G_2 \subset Z(G)$.  
\end{lemme}
\begin{proof} Notons tout d'abord que $V$ est nécessairement
irréductible, car sinon $V \otimes V$ aurait au moins
quatre composants. Supposons $G \simeq G_1 \times G_2$. Comme
$V$ est irréductible, $V \simeq V_1 \otimes V_2$, avec
$V_i$ une représentation irréductible de $G_i$. Comme $V$ est
fidèle, il en est de même des $V_i$. Or
$$
\begin{array}{lcl}
V \otimes V & \simeq & (V_1 \otimes V_2) \otimes (V_1 \otimes V_2) \\
& \simeq & (V_1 \otimes V_1) \otimes (V_2 \otimes V_2) \\
& \simeq & (S^2 V_1 \oplus  \Lambda^2 V_1) \otimes (S^2 V_2 \oplus  \Lambda^2 V_2) \\
& \simeq & (S^2 V_1 \otimes S^2 V_2) \oplus (S^2 V_1 \otimes \Lambda^2 V_2)
\oplus (\Lambda^2 V_1 \otimes S^2 V_2) \oplus (\Lambda^2 V_1 \otimes \Lambda^2 V_2) \\
\end{array}
$$
Donc, comme les dimensions des $V_i$ sont au moins 1,
on en déduit $\Lambda^2 V_2 = 0$ ou $\Lambda^2 V_1 = 0$, donc $V_2$ ou $V_1$
est de dimension 1. Comme les $V_i$ sont fidèles, on en déduit que
$G_1$ ou $G_2$ est commutatif, d'où la conclusion.
 
\end{proof}

On en déduit
\begin{proposition} \label{propdec} Soit $K$ un corps de caractéristique 0, et
$R : G \to GL_N(K)$ une re\-pré\-sen\-ta\-tion fidèle de $G$ telle que
$R(G)$ soit Zariski-dense dans $GL_N(K)$. Si $\Gamma$ est un sous-groupe
d'indice fini de $G$ et $\Gamma \simeq \Gamma_1 \times \Gamma_2$,
alors $\Gamma_1$ ou $\Gamma_2$ est commutatif.
\end{proposition}
\begin{proof} Soit $V = K^N$. On peut supposer que $K$ est clos parce que $R(G)$ est
également Zariski-dense dans $GL_N(\overline{K})$, $\overline{K}$ désignant
la clôture algébrique de $K$. Puisque l'adhérence de Zariski
$\overline{R(G)}$ de $R(G)$ est irréductible en tant que variété algébrique,
que $\Gamma$ est d'indice fini dans $G$ et que $R$ est fidèle, on a
$\overline{R(\Gamma)} = \overline{R(G)} = GL_N(K)$. Ainsi, $V$
et les deux composantes $S^2 V$ et $\Lambda^2 V$ de $V \otimes V$ sont
irréductibles sous l'action de $\Gamma$.
D'après le lemme précédent on en déduit la conclusion.
\end{proof}

Le théorème 5b est alors une conséquence immédiate du théorème 4 et de la proposition
précédente (voir également \cite{CORNHARP} pour une autre preuve de cette impliquation).
En fait, la même démonstration permet de montrer un résultat beaucoup plus fort,
puisque de nombreux sous-groupes de $B$ sont Zariski-denses dans un groupe linéaire.
Comme cette propriété de non-décomposabilité est une conséquence de la propriété (L)
du théorème 6, établie ci-dessous, nous nous contenterons de préciser les sous-groupes
pour lesquels nous avons démontré cette dernière propriété.

\subsection{Groupes fortement linéaires}
Si $K$ est un corps (commutatif) et $G$ un sous-groupe de $GL_N(K)$,
on désigne par $\overline{G}$ son adhérence dans $GL_N(K)$ pour la
topologie de Zariski. C'est un sous-groupe algébrique de $GL_N(K)$.
Rappelons que $Z(SL_N(K)) = \mu_N(K)$ est un groupe fini et que tout
sous-groupe (algébrique) distingué propre de $SL_N(K)$ est inclus dans
son centre. Suivant les conventions usuelles en théorie des groupes, $H < G$
(resp. $H \vartriangleleft G$) signifie que $H$ est un sous-groupe
(resp. un sous-groupe normal) de $G$ et, pour $H_1,H_2 < G$, la notation
$(H_1,H_2)$ désigne le sous-groupe de $G$ engendré par les commutateurs
$(x,y) = xyx^{-1} y^{-1}$ pour $x \in H_1, y \in H_2$. Lorsque
$H_1,H_2 \vartriangleleft G$, on a $H_1 \cap H_2 \subset (H_1,H_2)$.

\begin{definition} On dit qu'un groupe \emph{sans torsion} $G$ est \emph{fortement linéaire} en
dimension $N$ sur un corps (commutatif) $K$ s'il existe
un morphisme injectif $\iota : G \to GL_N(K)$
tel que $\overline{\iota(G)} \supset SL_N(K)$.
\end{definition}
S'il existe $N$ et $K$ tels que cette situation se présente, on dira
également que $G$ est fortement linéaire, sans préciser $N$ et $K$. Un
groupe fortement linéaire est évidemment linéaire au sens habituel.

Remarquons que, puisque $G$ est supposé sans torsion, tous ses
sous-groupes le sont é\-ga\-le\-ment. En particulier, identifiant un tel $G$
à son image dans $GL_N(K)$, pour $B < G$ la propriété $B < Z(SL_N(K))$
équivaut à $B = \{ e \}$. De plus, si $G \neq \{ e \}$ alors
$G$ est infini, donc $K$ est nécessairement un corps infini.

Si un tel $G$ est non abélien on a $N \geq 2$ et, inversement, si $G$
était abélien son action sur $K^N$ ne pourrait être irréductible pour
$N \geq 2$, ce qui contredirait $\overline{\iota(G)} \supset SL_N(K)$. Ainsi,
pour un tel groupe, $G$ abélien équivaut à $N=1$.

Rappelons les deux faits élémentaires suivants concernant la
topologie de Zariski.
\begin{enumerate}
\item si $K \subset L$ est une extension de corps, alors pour tout
$N \geq 1$ la topologie de Zariski de $K^N$ est induite par celle
de $L^N$
\item si de plus $K$ est infini, alors $K^N$ est Zariski-dense
dans $L^N$.
\end{enumerate}
En conséquence de (1), la topologie de Zariski de $GL_N(K)$ est induite
par celle de $GL_N(L)$. De plus, en tant que variété algébrique
sur $K$, $SL_N(K)$ est $K$-rationnelle, c'est-à-dire qu'elle admet un ouvert
isomorphe à un ouvert de l'espace affine. Une conséquence de (2) est donc
que, si $K$ est infini, alors $SL_N(K)$ est dense dans $SL_N(L)$.

Supposons que l'on ait un plongement $\iota_K : G \into GL_N(K)$, le groupe
$G$ étant toujours supposé sans torsion, donc en particulier $K$ est
un corps infini, et soit $L$ un surcorps de $K$.
On note $\iota_L : G \into GL_N(L)$ le morphisme induit.

\begin{lemme} \label{kclos} $\overline{\iota_K(G)} \supset SL_N(K)$ si
et seulement si
$\overline{\iota_L(G)} \supset SL_N(L)$.
\end{lemme}

\begin{proof} On peut supposer
$G < GL_N(K) < GL_N(L)$. On note pour les distinguer
$\overline{G}$ l'adhérence de $G$ dans $GL_N(K)$ et
$\widehat{G}$ son adhérence dans $GL_N(L)$. On a $\overline{G}
= \widehat{G} \cap GL_N(K)$ parce la topologie de Zariski de $GL_N(K)$
est induite de celle de $GL_N(L)$. Alors $\widehat{G} \supset SL_N(L)$
implique
$$
\overline{G} = \widehat{G} \cap GL_N(K) \supset
SL_N(L) \cap GL_N(K) = SL_N(K)
$$
et réciproquement $\overline{G} \supset SL_N(K)$ implique
$$
\widehat{G} = \widehat{\overline{G}} \supset \widehat{SL_N(K)} = SL_N(L).
$$
\end{proof}

En conséquence on pourra toujours supposer $K$ algébriquement
clos. De plus, sous ces hypothèses :
\begin{enumerate}
\item $G$ agit de façon absolument irréductible sur $K^n$.
\item $Z(G) \subset K^{\times}$, et $(G,G) \cap Z(G) = \{ e \}$.
\end{enumerate}

On rappelle d'autre part le résultat classique suivant. L'idée
de l'appliquer dans ce cadre est empruntée à \cite{CORNHARP}.

\begin{proposition} (cf. \cite{BOREL}, 2.4) \label{borel} Soit $\Gamma$ un
groupe algébrique défini sur
un corps $K$, et $M_1,M_2$ deux sous-groupes de $\Gamma(K)$
non nécessairement fermés
tels que $M_1$ normalise $M_2$. Alors $\overline{M_1}$ normalise $\overline{M_2}$ et
$\overline{(M_1,M_2)} = (\overline{M_1},\overline{M_2})$.
\end{proposition}

En particulier, un groupe fortement linéaire non abélien n'est jamais
nilpotent. En effet, notant $(C^r G)_{r \geq 0}$ la suite centrale
descendante du groupe $G$, une récurrence immédiate montre

\begin{ccor} Si $G$ est fortement linéaire en dimension $N$ sur $K$
il en est de même de $C^r G$ pour tout $r \geq 0$.
\end{ccor}

De nombreux sous-groupes des groupes fortement linéaires le sont
également.

\begin{proposition} \label{flinnorm} Soit $G$ un sous-groupe sans torsion
de $GL_N(K)$ tel que
$\overline{G} \supset SL_N(K)$. Pour tout $B \vartriangleleft G$ tel que
$B \not\subset Z(G)$ on a $\overline{B} \supset SL_N(K)$.
\end{proposition}

\begin{proof}
D'après le lemme \ref{kclos} on peut
supposer que $K$ est algébriquement clos. Comme $G$ normalise
$B$, d'après la proposition \ref{borel} on sait que
$\overline{G} \supset SL_N(K)$ normalise $\overline{B}$,
donc $\overline{B} \cap SL_N(K)$ est un sous-groupe
normal de $SL_N(K)$, algébrique sur $K$. 
Il suffit de montrer que $H = \overline{B} \cap SL_N(K) \not\subset \mu_N(K)$.
Or $(G,B) < B < \overline{B}$ car $B \vartriangleleft G$, de plus
$(G,B) < (GL_N(K),GL_N(K)) < SL_N(k)$ et enfin $(G,B) \neq \{ e \}$ d'après
l'hypothèse $B \not\subset Z(G)$. Ainsi $H$ contient le groupe
sans torsion $(G,B) \neq \{ e \}$ donc $H \not\subset \mu_N(K)$,
ce qui conclut.
\end{proof}

\begin{proposition} \label{flinfini}Soit $G$ un sous-groupe infini de $GL_N(K)$ tel que
$\overline{G} \supset SL_N(K)$. Pour tout $B < G$ d'indice
fini on a $\overline{B} \supset SL_N(K)$.
\end{proposition}
\begin{proof}
On suppose encore $K$ algébriquement clos. Soient $g_1,\dots,g_r \in G$
tels que $G = g_1 B \sqcup \dots \sqcup g_r B$. Alors
$\overline{G} = g_1 \overline{B} \cup \dots \cup g_r \overline{B}$ et
$SL_N(K) \subset \overline{G}$ implique que $SL_N(K)$ est la
réunion des $g_i \overline{B} \cap SL_N(K)$. Puisque $K$ est algébriquement
clos, c'est également vrai en tant que $K$-variété algébrique, donc il existe
$i \in [1,r]$ tel que $g_i \overline{B} \supset SL_N(K)$ par irréductibilité de
$SL_N$ comme variété algébrique. En particulier $e \in g_i \overline{B}$ donc
$g_i \in \overline{B}$ et $g_i \overline{B} = \overline{B} \supset SL_N(K)$.
\end{proof}

On rappelle qu'un sous-groupe $B$ d'un groupe $G$ est dit
\emph{sous-normal} s'il existe une famille $B_0,\dots B_n$ de
sous-groupes
de $G$ tels que $B = B_n < B_{n-1} < \dots < B_1 < B_0 = G$ et que,
pour tout $i \in [0,n-1]$, on ait $B_{i+1} \vartriangleleft B_i$.

\begin{theoC} Soit $G$ un groupe sans torsion fortement linéaire (en dimension $N$ sur $K$)
et $B$ un sous-groupe sous-normal de $G$ non inclus dans le centre de $G$. Alors
\begin{enumerate}
\item $B$ est fortement linéaire (en dimension $N$ sur $K$) ainsi que
ses sous-groupes d'indice fini.
\item $Z(B) = Z(G) \cap B$.
\item Si $A \vartriangleleft G$ et $A \not\subset Z(G)$, alors $A \cap B \not\subset Z(G)$.
\end{enumerate}
\end{theoC}
\begin{proof} On suppose $G < GL_N(K)$ avec $K$ algébriquement clos, et $\overline{G}
\supset SL_N(K)$. Si $G$ est abélien le résultat est évident donc on peut supposer
$N \geq 2$.
D'après le lemme de Schur on a $Z(G) \subset K^{\times}$. Soit $B$
un sous-groupe sous-normal de $G$ tel que $B \not\subset Z(G)$. On a donc 
$B = B_n < B_{n-1} < \dots < B_1 < B_0 = G$ avec $B_{i+1} \vartriangleleft B_i$,
et $B \not\subset Z(G)$. On démontrera, par récurrence sur $n$, les
assertions (2), (3) et (1'), avec (1') $\overline{B} \supset SL_N(K)$.
On a (1') $\Rightarrow$ (1) d'après la proposition \ref{flinfini}.  

Pour $n=0$ c'est immédiat. Supposons alors $B = B_n < G$ satisfaisant (1'), (2) et
(3), et soit $B' \vartriangleleft B$ tel que $B' \not\subset Z(G)$. On a $Z(B)
= Z(G) \cap B$ donc, puisque $B' \subset B$, on a $B' \not\subset Z(B)$. On déduit alors
de $\overline{B} \supset SL_N(K)$ que $\overline{B'} \supset SL_N(K)$ d'après la
proposition \ref{flinnorm}. Ainsi (1') est démontré. D'autre part (2) découle
de ce qu'alors $Z(B') \subset K^{\times}$, et $Z(G) = G \cap K^{\times}$
donc
$Z(G) \cap B' = G \cap K^{\times} \cap B' = K^{\times} \cap B' = Z(B')$.

Soit maintenant $A \vartriangleleft G$ tel que $A \not\subset Z(G)$, et $A' = A \cap B$.
Par hypothèse de récurrence on a $A' \not\subset Z(G)$
donc $A' \not\subset K^{\times}$ et en particulier $A' \not\subset Z(B)$.
Comme $A' \vartriangleleft B$ on peut appliquer la proposition \ref{flinnorm},
donc $\overline{A'} \supset SL_N(K)$. De plus $G$ donc $B$ et $B'$ normalisent
$A'$ d'où, d'après la proposition \ref{borel},
$$
\overline{(A',B')} = (\overline{A'},\overline{B'}) \supset (SL_N(K),
SL_N(K)) = SL_N(K).
$$
On en déduit $(A',B') \not\subset K^{\times}$ donc
$A' \cap B' \not\subset K^{\times}$ puisque $A',B' \vartriangleleft B$. Or
$Z(G) \subset K^{\times}$ donc $A' \cap B' \not\subset Z(G)$. On conclut
par récurrence.
\end{proof}

\begin{cor} Si $G$ est un groupe fortement linéaire contenant deux sous-groupes $A,B$
tels que $G \simeq A \times B$, alors $A \subset Z(G)$ ou $B \subset Z(G)$. Il en est
de même de ses sous-groupes sous-normaux et d'indice fini.
\end{cor}
\begin{proof} Par l'absurde, si $B \not\subset Z(G)$ et $A \not\subset Z(G)$
on peut appliquer le théorème puisque $A \vartriangleleft G$ et
$B \vartriangleleft G$. On en déduit $A \cap B \not\subset Z(G)$, ce qui
est absurde puisque $A \cap B = \{ e \}$.
\end{proof}

Notons $\mathrm{Fit}(G)$ le sous-groupe de Fitting de $G$, c'est-à-dire
le sous-groupe engendré par les sous-groupes normaux nilpotents de $G$.
  
\begin{cor} Si $G$ est un groupe fortement linéaire, alors $\mathrm{Fit}(G)
=Z(G)$.
\end{cor}
\begin{proof} En effet, $Z(G) < \mathrm{Fit}(G)$ et, si $N \vartriangleleft
G$ est nilpotent, alors $N$ ne peut être fortement linéaire d'après
le corollaire de la proposition \ref{borel}, donc $N \subset Z(G)$
et $\mathrm{Fit}(G) = Z(G)$.
\end{proof}

Notons $\Phi(G)$ le sous-groupe de Frattini de $G$, intersection de ses
sous-groupes maximaux.

\begin{cor} Si $G$ est un groupe fortement linéaire qui est finiment
engendré, alors $\Phi(G) \subset Z(G)$.
\end{cor}

\begin{proof} Le cas où $G$ est abélien est évident, donc on peut
l'exclure. Comme $G$ est linéaire et finiment engendré, $\Phi(G)$
est nilpotent d'après le théorème de Platonov \cite{PLATO}. Comme
$\Phi(G) \vartriangleleft G$ on en déduit $\Phi(G) \subset Z(G)$,
car sinon $\Phi(G)$ serait fortement linéaire en dimension $N \geq 2$,
contredisant sa nilpotence.
\end{proof}

\section{Etude du type $A_{n-1}$}

\subsection{Généralités}

Dans le cas où $W \simeq \SN$ est de type $A_{n-1}$, on choisit pour $\widetilde{W} \subset W$
un parabolique naturel de type $A_{n-2}$, isomorphe à $\mathfrak{S}_{n-1}$.
On note $\mathcal{T}_n = \g$, et $\mathcal{T}_{n-1} = \widetilde{\g}$
l'algèbre de Lie associée à $\widetilde{W}$.

La représentation $V_n = V$ de $\mathcal{T}_n$ a pour dimension
$n(n-1)/2$. Notons $v_{ij}$
et $t_{ij}$ respectivement le vecteur de base et le générateur
correspondant à la transposition $(i \ j)$. L'action est alors décrite
par les formules
$$
\left\lbrace \begin{array}{lcl}
t_{ij}.v_{jk} & = & v_{ik} - v_{ij} \\
t_{ij}.v_{kl} & = & v_{kl} \\
t_{ij} . v_{ij} & = & m v_{ij} 
\end{array} \right.
$$
où $i,j,k,l$ correspondent à des indices distincts.
On note $V_n$ l'espace vectoriel sous-jacent, engendré par les $v_{ij}$,
et $\mathcal{T}_n$ l'algèbre de Lie des tresses infinitésimales pures,
engendrée par les $t_{ij}$. On note $T_n = \sum t_{ij}$, qui agit
par $m + c(W)$ sur $V_n$, et on plonge de façon naturelle $\mathcal{T}_{n-1}$
dans $\mathcal{T}_n$.

\subsection{Restriction à $\mathcal{T}_{n-1}$}

Le $\mathcal{T}_{n-1}$-module $V_n$ admet pour sous-espace stable
le sous-espace engendré par les $v_{ij}$ pour $i,j < n$. En
tant que $\mathcal{T}_{n-1}$-module, il est naturellement
isomorphe à $V_{n-1}$. 

Notons $T' = T_{n-1}$, et, pour $ 1 \leq k \leq n-1$,
$$
v'_k = \sum_{\stackrel{i \neq k}{i \leq n-1}} v_{ik}, \ \ \ 
w^0_{kn} = v_{kn} + \frac{1}{m-n+4} v'_k,
w^0 = \sum_{i \leq n-1} w^0_{in}, \check{v} = \sum_{i,j \leq n-1} v_{ij}.
$$
On introduit alors les sous-espaces
suivants

$$
\begin{array}{|lcl|c|}
\hline
\multicolumn{3}{|c|}{\mathrm{Sous-espaces}} & \mathrm{Dimension} \\
\hline
V_{n-1} & = & < v_{ij} \ | \ 1 \leq i,j \leq n-1 > & \# \widetilde{\RR} \\
\hline
U_{n-1} & = & < w^0_{kn} - w^0_{k+1,n} \ | \ 1 \leq k \leq n-2 > & n-2 \\
\hline
S_{n-1} & = & \displaystyle < w^0 + \frac{2(n-1)}{(m-n+4)(m-2n+5)} \check{v} > & 1 \\
\hline
\end{array}
$$

Chacun de ces sous-espaces est stable pour
l'action de $\mathcal{T}_{n-1}$, et on a $V_n = V_{n-1} \oplus U_{n-1}
\oplus S_{n-1}$.
Une base commode de $U_{n-1} \oplus S_{n-1}$ est donnée par les
vecteurs $w^1_{kn} = w^0_{kn} + \frac{2}{(m-n+4)(m-2n+5)} \check{v}$
pour $1 \leq k \leq n-1$. Dans cette base, on a 
$
t_{ij}.w^1_{kn} = w^1_{kn}$ si $k \not\in \{i,j \}$,
$t_{ik} . w^1_{kn} = w^1_{in}$. On en déduit que $U_{n-1} \oplus S_{n-1}$
s'identifie à la représentation de Burau (non réduite) infinitésimale
de \cite{QUOT,IRRED}.

\subsection{Forme quadratique}

On suppose $n \geq 3$.
On a $(v_{ij}| v_{ij}) = m-1$,
$(v_{ij}| v_{jk}) = -1$ si $\#\{i,j,k\} = 3$,
$(v_{ij}| v_{k,l}) = 0$ si $\#\{i,j,k,l\} = 4$. Soient
$1 \leq i,j ,k\leq n-1$ et $i \neq j$. On a
$$
(v_{ij}| v'_k) = \left\lbrace
\begin{array}{ll} m-n+2 & \mbox{si } k \in \{i,j \} \\
-2 & \mbox{si } k \not\in \{i,j \} 
\end{array} \right. \ \   
(v_{ij}| v_{kn}) = \left\lbrace
\begin{array}{ll} -1 & \mbox{si } k \in \{i,j \} \\
0 & \mbox{si } k \not\in \{i,j \} 
\end{array} \right. \ \   
(v_{ij}|\check{v}) = m-2n +5
$$
On en déduit
$$
(w^0_{kn}| v_{ij}) = \frac{-2}{m-n+4}, (w^0_{kn}| \check{v})
= \frac{-(n-1)(n-2)}{m-n+4}, (\check{v}|\check{v}) =
(m-2n+5) \frac{(n-1)(n-2)}{2}
$$
et en particulier que le sous-espace $V_{n-1}$ est orthogonal à
$U_{n-1}\oplus S_{n-1}$. D'autre part, si
$1 \leq k \leq n-1$,
$$
\begin{array}{|c||c|c|c|}
\hline
 & (v_{kn}| v_{ln}) & (v_{kn}| v'_{l}) & (v'_{k}| v'_{l}) \\
\hline
\hline
k=l & m-1 &2-n & (n-2)(m-n+2) \\
\hline
k \neq l & -1 & -1 & m-3n+8 \\
\hline \end{array}
$$
On en déduit $(w^1_{kn}| w^1_{ln}) = \alpha$ si $k=l$, $\beta$
si $k \neq l$, avec
$$
\alpha = 
\frac{(m^2+(5-2n)m-2)(m-n+3)}{(m-2n+5)(m-n+4)},
\beta = 
\frac{-(m-2n+7)(m-n+3)}{(m-2n+5)(m-n+4)}
$$
Ainsi la forme quadratique $Q_n$ sur $V_n$ est
somme directe de $Q_{n-1}$ et d'une forme
quadratique $\delta_{n-1}$, que l'on peut
définir pour $m \not\in \{ n-4,2n-5 \}$, de discriminant
$(\alpha - \beta)^{n-2}(\alpha + (n-2)\beta)$  soit
$$
\frac{(m+1)^{n-2}(m-n+3)^{n-1}(m-2n+3)}{(m-n+4)^{n-2}(m-2n+5)}
$$
Comme $Q_n = Q_{n-1} \oplus \delta_{n-1}$,
et que $Q_2$ est la forme quadratique sur
la droite engendrée par $v_{12}$ de discriminant
$m-1$, on en déduit par une récurrence facile
\begin{proposition} \label{discrA} Pour $W$ de type $A_{n-1}$ avec $n \geq 3$, le discriminant de la forme quadratique
$(\ | \ )$
est
$$
(m+1)^{\frac{n(n-3)}{2}} (m-n+3)^{n-1}(m-2n+3) $$
et cette forme quadratique est définie positive ssi $m > 2n-3$.
\end{proposition}
La condition de définie positivité provient de la remarque générale
sur $(\ | \ )$ lors de sa définition (section 3.2).

\section{Etude du type $D_n$}

\subsection{Généralités} On suppose $n \geq 4$, et que $W$ est de type
$D_n$. 
Notons $(i \ j)$ et $(i\ j)'$ les réflexions respectivement définies par
$$
\begin{array}{lclcl}
(i\ j) & : & (z_1,\dots,z_i,\dots, z_j,\dots,z_n) &
\mapsto & (z_1,\dots,z_j,\dots, z_i,\dots,z_n) \\
(i\ j)'& : & (z_1,\dots,z_i,\dots, z_j,\dots,z_n)
& \mapsto & (z_1,\dots,-z_j,\dots,-z_i,\dots,z_n)
\end{array}
$$
$v_{ij}$ et
$v'_{ij}$, $t_{ij}$, $t'_{ij}$ les vecteurs de base et les générateurs
correspondant respectivement à
$(i \ j)$ et $(i \ j)'$. L'action est alors décrite par les formules
$$
\left\lbrace \begin{array}{lcl}
t_{ij} . v_{kl} & = & v_{kl} \\
t_{ij} . v_{jk} & = & v_{ik} - v_{ij} \\
t_{ij} . v_{ij} & = & mv_{ij} \\
t_{ij} . v'_{kl} & = & v'_{kl} \\
t_{ij} . v'_{jk} & = & v'_{ik}- v_{ij} \\
t_{ij} . v'_{ij} & = & v'_{ij} \\
\end{array} \right.
\left\lbrace \begin{array}{lcl}
t'_{ij} . v_{kl} & = & v_{kl} \\
t'_{ij} . v_{jk} & = & v'_{ik} - v'_{ij} \\
t'_{ij} . v_{ij} & = & v_{ij} \\
t'_{ij} . v'_{kl} & = & v'_{kl} \\
t'_{ij} . v'_{jk} & = & v_{ik}- v'_{ij} \\
t'_{ij} . v'_{ij} & = & mv'_{ij} \\
\end{array} \right.
$$
où $i,j,k,l$ désignent des indices distincts. On note $V_n = V$
l'espace vectoriel sous-jacent engendré par les $v_{ij}$ et $v'_{ij}$,
et $\mathcal{T}^D_n = \g$ est
engendrée par les $t_{ij}$ et
$t'_{ij}$. On considère $\mathcal{T}^D_{n-1} \subset \mathcal{T}^D_n$,
correspondant au sous-groupe parabolique maximal $\widetilde{W}$
engendré par les $(i\ j)$, $(i\ j)'$ pour $1 \leq i, j \leq n-1$, et
$V_{n-1} \subset V_n$ la représentation associée.

\subsection{Restriction à $\mathcal{T}^D_{n-1}$}

On définit, pour tout
$1 \leq k \leq n-1$,
$$
u_k = \sum_{\stackrel{i \leq n-1}{i \neq k}} v_{ik} + v'_{ik},
u = \sum_{\{r,s \} \subset [1,n-1]} v_{rs} + v'_{rs}, w_k = v'_{kn} - v_{kn}
$$
et $q_k = (4n-m-11)(2n-m-9)(v_{kn} + v'_{kn}) -2(4n-m-11) u_k + 8 u$.
Alors
$$
\left\lbrace \begin{array}{rclcrcl}
t_{ij} . q_k & = & t'_{ij} . q_k = q_k  \\
t_{ij} . q_j & = & t'_{ij} . q_j = q_i 
\end{array} \right.
\left\lbrace \begin{array}{rclcrcl}
t_{ij} . w_k & = & w_k ,&  \ & t'_{ij} . w_k & = & w_k \\
t_{ij} . w_j & = & w_i ,&  \ & t'_{ij} . w_j & = & -w_i 
\end{array} \right.
$$
où $i,j,k$ sont distincts et compris entre $1$ et $n-1$. Si
$m \not\in \{ 4n-11,2n-9 \}$, on a une décomposition en somme directe
de $\mathcal{T}_{n-1}^D$-modules, irréductibles pour $m$ générique, $V_n = V_{n-1} \oplus U^D_{n-1} \oplus
U^A_{n-1} \oplus S_{n-1}$ avec 
$$
\begin{array}{|lcl|c|}
\hline
\multicolumn{3}{|c|}{\mathrm{Sous-espaces}} & \mathrm{Dimension} \\
\hline
V_{n-1} & = & < v_{ij}, v'_{ij} \ | \ 1 \leq i,j \leq n-1 > & \# \widetilde{\RR} \\
\hline
U_{n-1}^D & = & < w_k \ | \ 1 \leq k \leq n-1 > & n-1 \\
\hline
U_{n-1}^A & = & < q_k - q_{k+1} \ | \ 1 \leq k \leq n-2 > & n-2 \\
\hline
S_{n-1} & = & < q_1 + \dots + q_{n-1} > & 1 \\
\hline
\end{array}
$$

On remarque que l'espace $U^A_{n-1} \oplus S$ est l'intersection des noyaux
des $t_{ij} - t'_{ij}$, pour $1 \leq i,j \leq n-1$. D'autre part
$U_{n-1}^D$ est la représentation de Burau infinitésimale de type $D_{n-1}$,
et $U_{n-1}^A$ est la représentation de Burau infinitésimale de type
$A_{n-2}$, étendue à $\mathcal{T}^D_{n-1}$ par le morphisme
d'algèbres de Lie $\mathcal{T}^D_{n-1} \to \mathcal{T}_{n-1}$
défini par $t_{ij},t'_{ij} \mapsto t_{ij}$.

\subsection{Forme quadratique}

La forme quadratique est définie par
$(v_{ij}|v_{ij}) = (v'_{ij}|v'_{ij}) = m-1$, $(v_{ij}|v_{kl})
= (v_{ij}|v'_{kl}) = (v'_{ij}|v'_{kl}) = 0$, 
$(v_{ij}|v'_{ij}) = 0$ et $(v_{ij}|v_{jk})
= (v_{ij}|v'_{jk}) = (v'_{ij}|v'_{jk}) = -1$ si $\#\{ i,j,k,l\} = 4$.

On constate facilement que $(w_k|w_k) = 2(m-1)$ et $(w_i|w_j)=0$ si
$i \neq j$. On en déduit que le discriminant de $(\ | \ )$ sur
$U^D_{n-1}$ est $2^{n-1}(m-1)^{n-1}$. La restriction de $( \ | \ )$
sur $U^A_{n-1}\oplus S$ est donnée par $( q_k | q_k ) = \alpha$,
$( q_i | q_j ) = \beta$ si $i \neq j$, avec
$$
\begin{array}{l}
\alpha = -2(m-2n+9)(4n+4nm-12m-3-m^2)(m-4n+11)(m-2n+7) \\
\beta = -4(m-2n+9)(m-2n+7)(m-4n+15)(m-4n+11)
\end{array}
$$
donc son discriminant est $(\alpha-\beta)^{n-2}(\alpha + (n-2) \beta)$,
avec
$$
\begin{array}{l}
\alpha - \beta =  
2(m+3)(m-2n+7)(m-2n+9)(m-4n+11)^2 \\
\alpha + (n-2) \beta =
2(m-4n+7)(m-2n+7)(m-4n+11)(m-2n+9)^2
\end{array}
$$

On en déduit comme pour le type A, par une récurrence facile, que

\begin{proposition} \label{discrD} En type $D_n$, $n \geq 3$, le discriminant de la
forme quadratique $(\ | \ )$ vaut
$$ (m-4n+7)(m-1)^{\frac{n(n-1)}{2}} (m+3)^{\frac{n(n-3)}{2}} (m-2n+7)^{n-1}
$$
et elle est définie positive ssi $m > 4n-7$.
\end{proposition}

\end{document}